\documentclass{article}

\usepackage{arxiv}

\usepackage[utf8]{inputenc} % allow utf-8 input
\usepackage[T1]{fontenc}    % use 8-bit T1 fonts
\usepackage{hyperref}       % hyperlinks
\usepackage{url}            % simple URL typesetting
\usepackage{booktabs}       % professional-quality tables
\usepackage{amsfonts}       % blackboard math symbols
\usepackage{nicefrac}       % compact symbols for 1/2, etc.
\usepackage{microtype}      % microtypography
\usepackage{lipsum}		% Can be removed after putting your text content
\usepackage{graphicx}
\usepackage[square,numbers]{natbib}
\usepackage{doi}
\usepackage{amsmath}
\usepackage{siunitx}
\usepackage{amssymb}
\usepackage{graphicx}
\usepackage{subcaption}
\usepackage{bm} % bold math symbols
\usepackage{xspace}
\usepackage{mathtools}
\usepackage{tikz, smartdiagram}
\usepackage[nameinlink]{cleveref} % for smart reference including capitalization

\newcommand{\velocity}{\textbf{v}}
\newcommand{\ddx}[2]{\frac{\partial #1}{\partial #2}}

\newcommand\norm[1]{\left\lVert#1\right\rVert}

\DeclareMathOperator*{\argmin}{arg\,min}

\newcommand{\fenics}{\texttt{FEniCs}\xspace}
\newcommand{\hippylib}{\texttt{hIPPYlib}\xspace}

\newcommand{\inflowBoundary}{\Gamma_{-}}
\newcommand{\characteristicBoundary}{\Gamma_{0}}
\newcommand{\outflowBoundary}{\Gamma_{+}}

\newcommand{\diffcoeff}{\kappa}

\newcommand{\parameterfield}{m}

\newcommand{\normal}{\textbf{n}}

\newcommand{\spacepoint}{\textbf{x}}

\newcommand{\stdbase}{\textbf{e}}
\newcommand{\measurement}{\mathbf{d}}

\newcommand{\obsO}{\mathcal{B}} % observation operator 
\newcommand{\pto}{\mathcal{F}} % parameter to observable map
\newcommand{\ptoh}{\pto_h} % parameter to observable map discrete 
\newcommand{\pts}{\mathcal{K}} % parameter to state operator
\newcommand{\sop}{D} % space of parameter (FE)
\newcommand{\ndof}{{n_{\text{dof}}}}
\newcommand{\sensorweights}{\mathbf{w}}
\newcommand{\W}{W}
\newcommand{\Wsqrt}{\W^{\frac{1}{2}}}
\newcommand{\hessian}{\mathcal{H}}
\newcommand{\hessianh}{\hessian_h}
\newcommand{\hessianhmisfit}{\hessian^{\textbf{misfit}}_h}
\newcommand{\covPr}{\Gamma_{\text{pr}}}
\newcommand{\covPrh}{\Gamma_{\text{pr,h}}}
\newcommand{\precondhessianhmisfit}{\tilde\hessian^{\textbf{misfit}}_h}

\newcommand{\QoiO}{\mathcal{P}}
\newcommand{\Tzeroqoi}{T^{\text{QoI}}_0}
\newcommand{\Tfinalqoi}{T^{\text{QoI}}_\text{final}}
\newcommand{\Qoi}{P}

\newcommand{\rpto}{r_\pts}% forward equation operator
\newcommand{\R}{\mathbb{R}}

\newcommand{\tobs}{t^{\text{obs}}}
\newcommand{\xobs}{\textbf{x}^{\text{obs}}}

\newcommand{\observationforequation}{({\tobs_i,\spacepoint_i^\text{obs}})}

\newcommand{\ansatzSpace}{\mathcal{V}_h}

\newcommand{\FEMi}{I}
\newcommand{\ellipI}{\eta}
\newcommand{\ellipLaplace}{\gamma}

\newcommand{\fracSolidus}[2]{#1/#2\,}

\newcommand{\gauss}[2]{\mathcal{N}(#1,#2)}
\newfont{\amsbold}{msbm10}
\newfont{\logobold}{logobf10 scaled\magstep2}

\def\*#1{\mathbf{#1}}

\newcommand{\bd}{\mathbf{d}}

\newcommand{\bw}{\mathbf{w}}

\newcommand{\by}{\mathbf{y}}

\newcommand{\mmap}{m_{\text{map}}}
\newcommand{\Tstep}{T_{\text{step}}}

\newcommand{\indFunc}[1]{\displaystyle \mathbf {1}_{#1}}

\title{Goal-oriented optimal sensor placement for PDE-constrained inverse problems in crisis management}

\date{} 					% Or removing it

\author{ \href{https://orcid.org/0009-0009-6325-3578}{\includegraphics[scale=0.06]{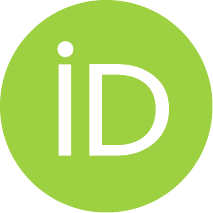}\hspace{1mm}Marco Mattuschka}\thanks{Corresponding author}\\
    German Aerospace Center (DLR)\\
	Institute for the Protection of Terrestrial Infrastructures\\
	53757 Sankt Augustin, Germany \\
	\texttt{marco.mattuschka@dlr.de} \\
	\And
    Noah An der Lan \\
    German Aerospace Center (DLR)\\
	Institute for the Protection of Terrestrial Infrastructures\\
	53757 Sankt Augustin, Germany \\
    \And
    \href{https://orcid.org/0000-0002-2814-0027}{\includegraphics[scale=0.06]{orcid.pdf}\hspace{1mm}Max von Danwitz} \\
    German Aerospace Center (DLR)\\
	Institute for the Protection of Terrestrial Infrastructures\\
	53757 Sankt Augustin, Germany \\
    \And
    \href{https://orcid.org/0000-0001-5767-7803}{\includegraphics[scale=0.06]{orcid.pdf}\hspace{1mm}Daniel Wolff} \\
    University of the Bundeswehr Munich\\
	Institute for Mathematics and \\Computer-Based Simulation (IMCS)\\
	85577 Neubiberg, Germany\\
    \And
    \href{https://orcid.org/0000-0002-8820-466X}{\includegraphics[scale=0.06]{orcid.pdf}\hspace{1mm}Alexander Popp} \\
    University of the Bundeswehr Munich\\
	Institute for Mathematics and \\Computer-Based Simulation (IMCS)\\
	85577 Neubiberg, Germany\\
    German Aerospace Center (DLR)\\
	Institute for the Protection of Terrestrial Infrastructures\\
	53757 Sankt Augustin, Germany \\
}
\renewcommand{\headeright}{}
\renewcommand{\undertitle}{}
\renewcommand{\shorttitle}{Goal-oriented optimal sensor placement}

\hypersetup{
pdftitle={A template for the arxiv style},
pdfsubject={q-bio.NC, q-bio.QM},
pdfauthor={David S.~Hippocampus, Elias D.~Striatum},
pdfkeywords={First keyword, Second keyword, More},
}

\begin{document}
\maketitle

\begin{abstract}
This paper presents a novel framework for goal-oriented optimal static sensor placement and dynamic sensor steering in PDE-constrained inverse problems, utilizing a Bayesian approach accelerated by low-rank approximations. The framework is applied to airborne contaminant tracking, extending recent dynamic sensor steering methods to complex geometries for computational efficiency. A C-optimal design criterion is employed to strategically place sensors, minimizing uncertainty in predictions. Numerical experiments validate the approach’s effectiveness for source identification and monitoring, highlighting its potential for real-time decision-making in crisis management scenarios.
\end{abstract}

% keywords can be removed
\keywords{Airborne contaminant transport \and Large-scale inverse problems \and Optimal experimental design \and Dynamic sensor steering}

\section{Introduction}

\begin{figure}
\centering
\includegraphics[width=0.45\linewidth]{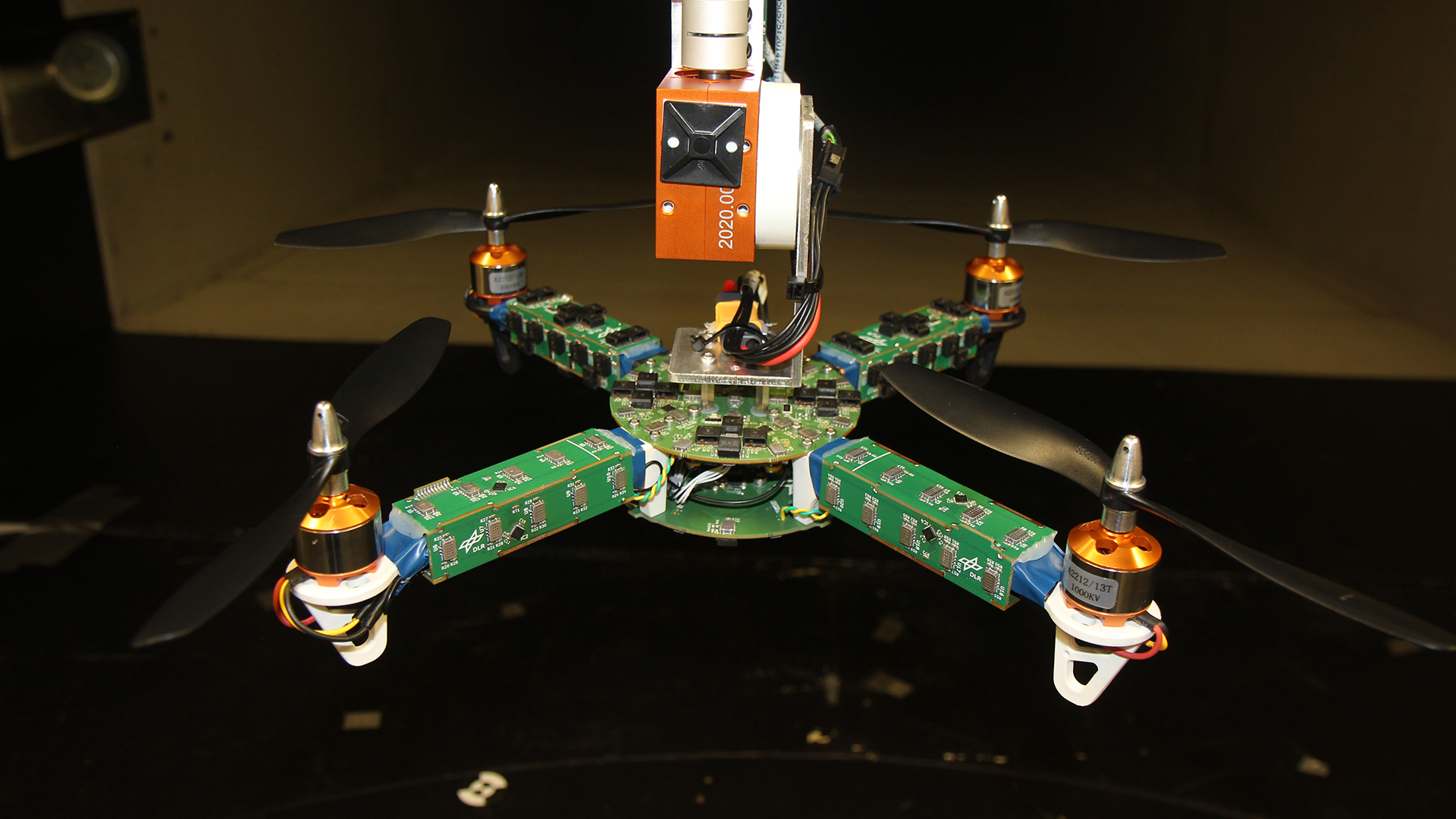} \hspace{0.5cm}
\includegraphics[width=0.45\linewidth]{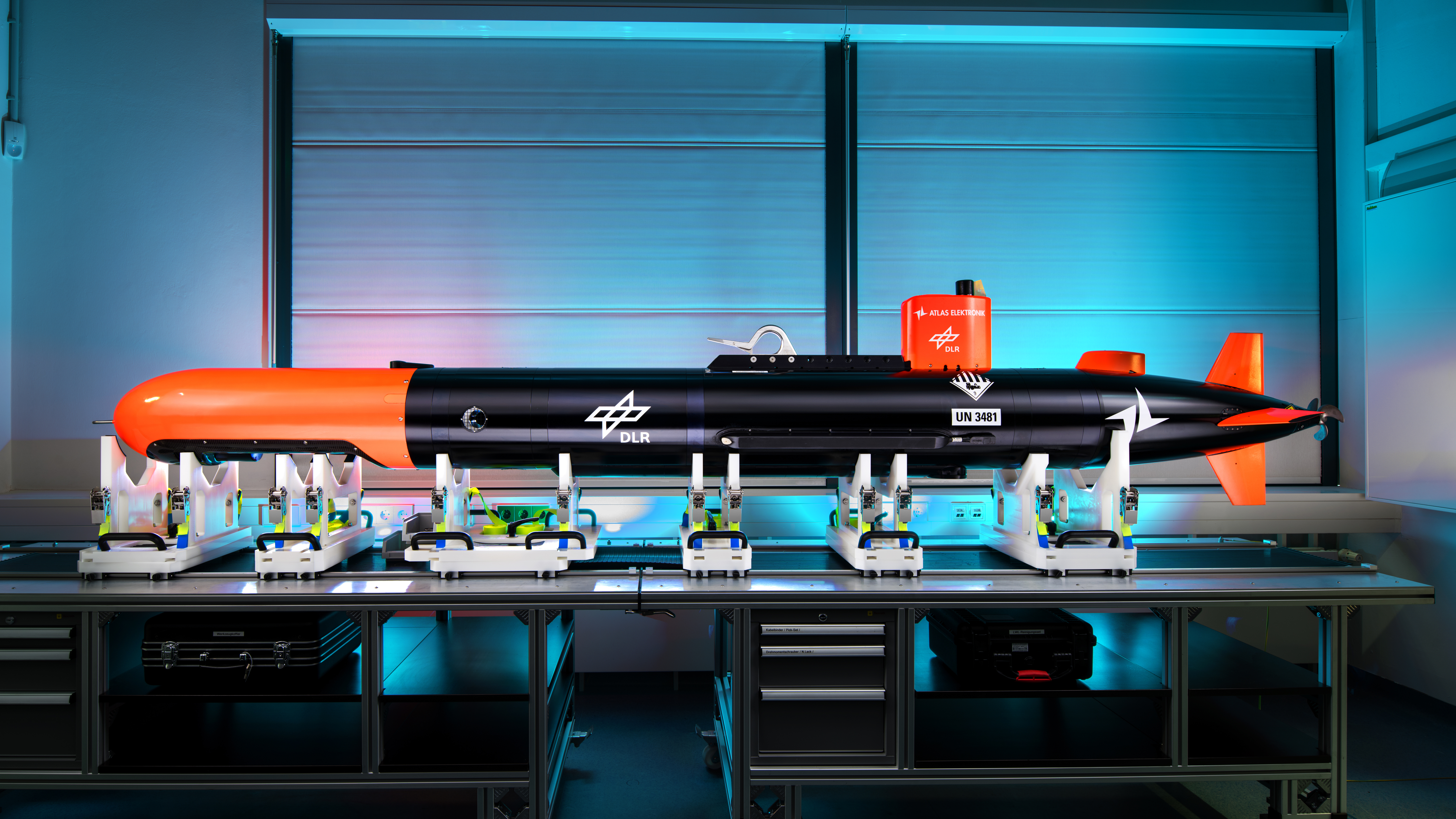} \\ \vspace{0.5cm}
\includegraphics[width=0.45\linewidth]{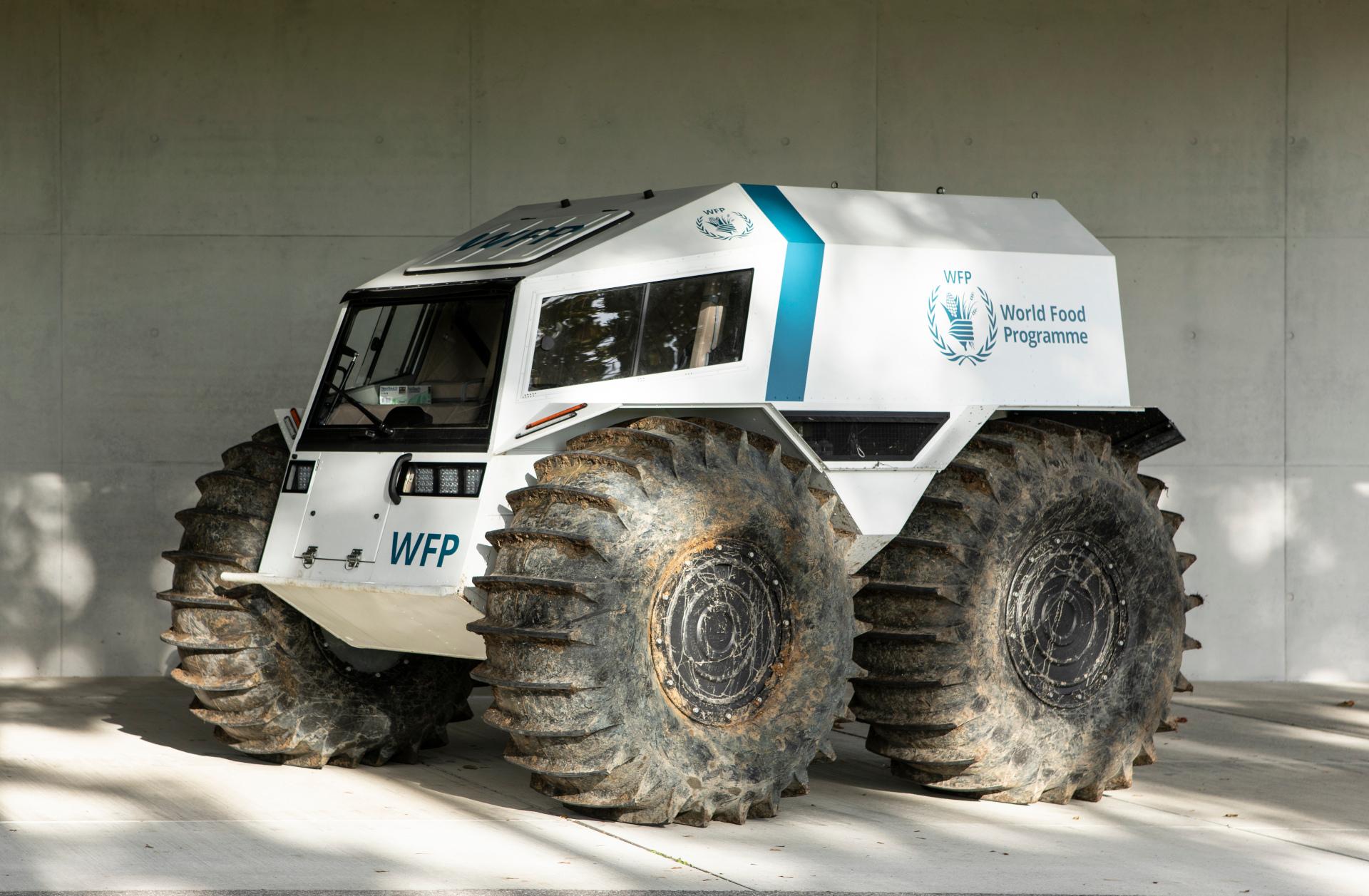} \hspace{0.5cm}
\includegraphics[width=0.45\linewidth]{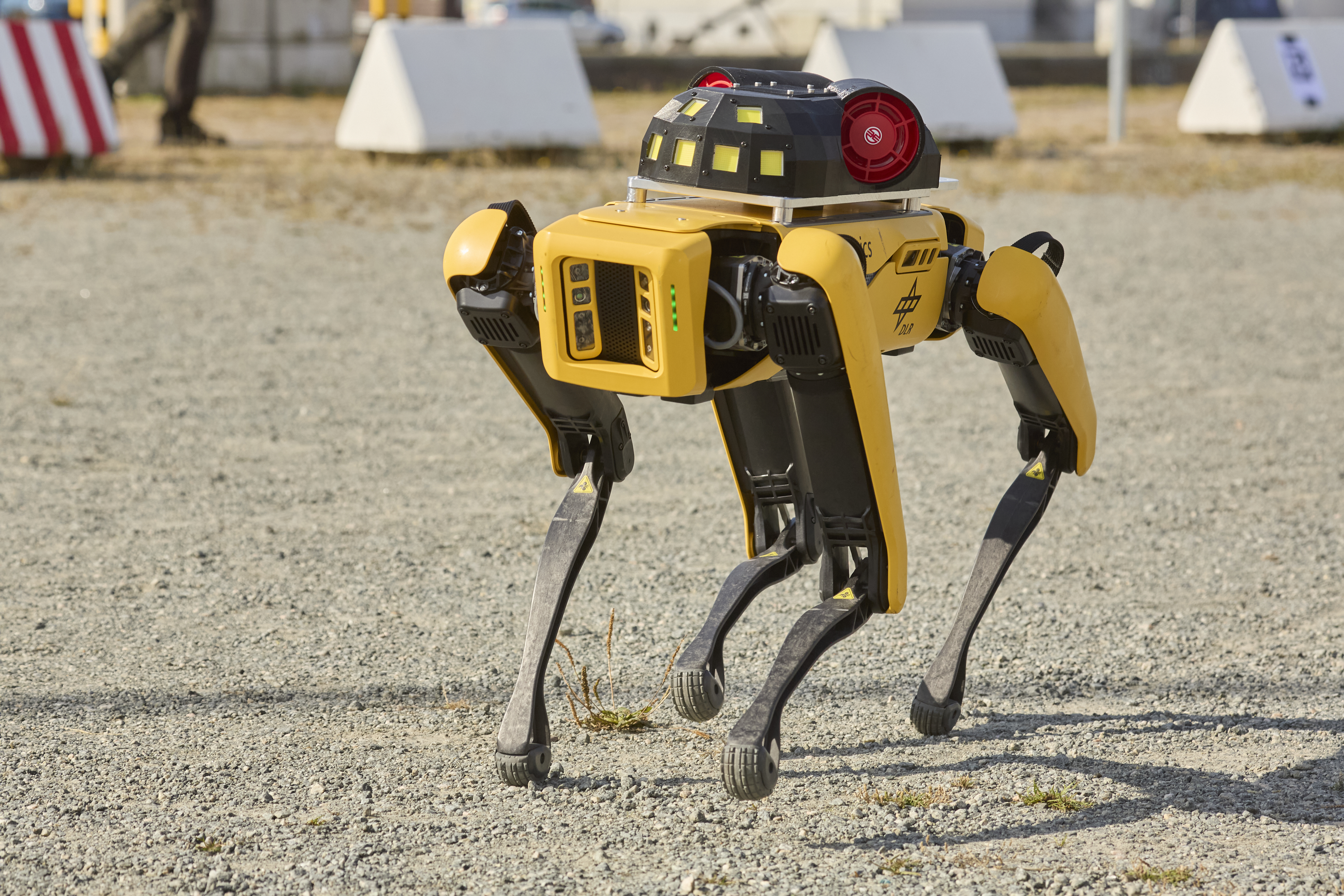}
\caption{Unmanned systems operated as autonomous sensor platforms at German Aerospace Center (DLR). Images: DLR, CC BY-NC-ND 3.0}
\label{fig:uxs}
\end{figure}

The growing capabilities of unmanned systems (UxS), such as unmanned aerial vehicles (UAVs, drones), autonomous underwater vehicles (AUVs), and unmanned ground vehicles, e.g., mobile robots, have made these systems indispensable tools in crisis management situations (see \autoref{fig:uxs} for examples). Such autonomous sensor platforms facilitate the collection of valuable information in regions where manned missions would be too dangerous or simply impossible due to inaccessibility for humans. However, during the routing of a sensor platform in the complex environment of a crisis situation, the question arises as to at which locations (additional) measurements provide an information gain and thereby constitute added value for decision makers. This very question leads to the need for a goal-oriented optimal sensor placement and in the dynamic case to a sensor steering problem.

This contribution addresses the problem in the mathematical setting of inverse problems constrained by partial differential equations (PDEs). A selected sensor steering strategy is developed and applied to the specific challenge of mapping airborne contaminant dispersion in the region of interest using discrete sensor measurements. Current methods for contaminant source identification and spread prediction rely heavily on an informative sensor placement. The selection of measuring points (sensor positions) is crucial, yet many existing studies focus solely on stationary sensors. This work bridges this gap by incorporating recent advances in sensor selection and experimental design to derive an algorithm for optimal sensor steering. Our research focuses on developing a systematic approach to select sensor positions that maximize the accuracy of contaminant source identification and prediction. By integrating current methodological advancements, our aim is to provide a practical solution for rescue workers and first responders, allowing informed decisions in high-stakes situations. Whereas the numerical examples presented in this work focus on the specific application of airborne contaminant transport, the goal-oriented optimal sensor placement strategy is independent of the considered physical model and, hence, can be easily transferred to other crisis management applications.

Starting point for this work is the recent publication by Wogrin et al.\cite{Wogrin.2023} that pioneers a dynamic sensor steering method in the context of airborne contaminant transport. In this work, we extend this approach to a significantly more complex geometry. Moreover, the inverse problem solution follows a more advanced approach that uses a Laplacian-like operator of trace class as prior information within a Bayesian inverse problem framework \cite{Villa.2021, Petra.2011}. To achieve approximate real-time capability, low-rank approximations of the Hessian matrix are precomputed in an offline phase \cite{Halko.2011, Liberty.2007}, enabling efficient problem solving in the online phase using a preconditioned inexact Newton-CG solver \cite{Steihaug.1983}. 

In the calculation of an optimal (stationary) sensor layout, Alexanderian et al.~\cite{Alexanderian.2018} used a reduced model for the contaminant transport to determine an A-optimal design that minimizes the average point-wise posterior variance of the inferred parameter vector. Following extensions proposed in \cite{Spantini.2017} and \cite{Attia.2018}, we focus on a goal-oriented design, i.e., the uncertainty associated with the prediction of the contaminant concentration in a specific region and time is minimized. To demonstrate the feasibility of sensor steering, a relaxed optimality criterion compared to the A-optimal criterion is chosen for this test case. Specifically, we use the C-optimal criterion, which focuses on minimizing the posterior variance of a particular linear combination of the inversion parameters. This approach eliminates the need to estimate the trace of the full covariance matrix and allows us to directly assess the impact of the covariance matrix on the parameter of interest. Alternatively, a D-optimal goal-oriented design in infinite dimensions maximizes the expected information gain~\cite{Wu.2023, Alexanderian.2018}. For a broader perspective on optimal experimental design (OED) for infinite-dimensional Bayesian inverse problems governed by PDEs, the interested reader is referred to \cite{Alexanderian.2021}.
The remainder of this paper is organized as follows. \cref{sec:background} provides background and mathematical formulations of the forward problem of contaminant transport, the inverse problem of source identification, as well as sensor positioning strategies and goal-oriented optimization. The combination of methodological developments into an algorithm for goal-oriented optimal sensor placement and steering is described in \cref{sec:method}. Numerical results are presented in \cref{sec:numerical} for three test cases of goal-oriented optimal sensor placement, namely (a) to identify an instantaneous contaminant source in a user-defined area of interest, (b) to monitor an area of special interest over a predefined time window, and (c) to steer a moving sensor. Finally, \cref{sec:conc} offers a conclusion and an outlook.

\section{Background}\label{sec:background}

\subsection{Forward problem: Contaminant distribution evolution} \label{ssec:ad-forward}
\begin{figure}[ht!]
\includegraphics[width=0.3\linewidth]{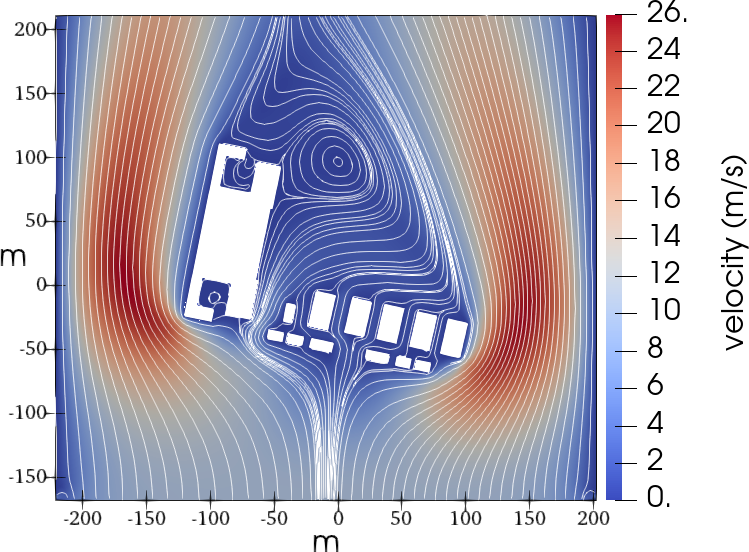}
\includegraphics[width=0.33\linewidth]{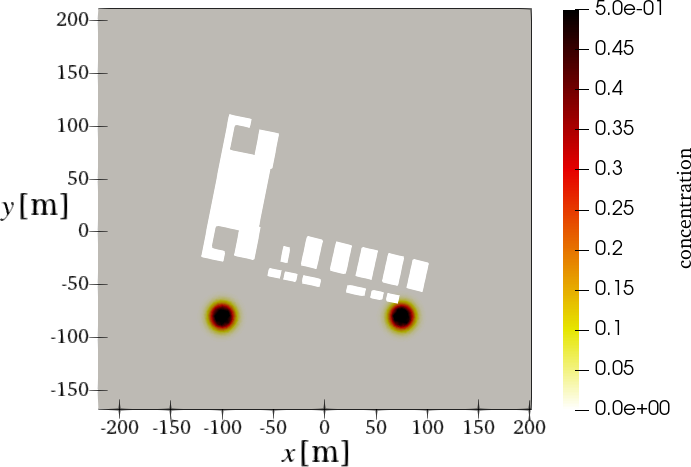}
\includegraphics[width=0.33\linewidth]{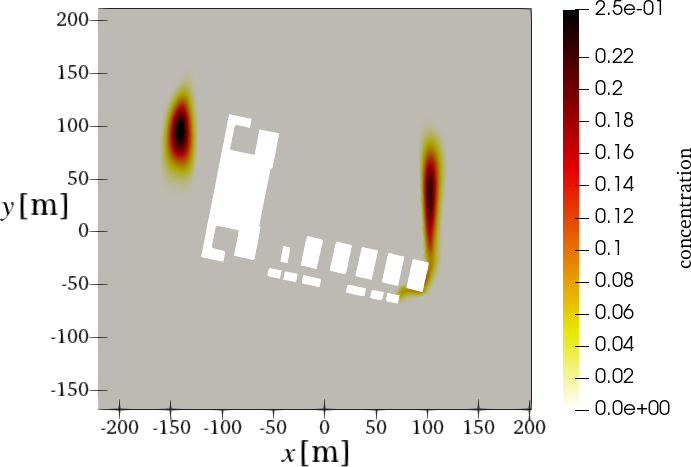}
\caption{Forward simulation of airborne contaminant transport on a campus geometry. Estimated wind vector field $\velocity$ (left), initial condition (middle), and simulated concentration at $t= 10\,\si{\second}$ (right).}
\label{fig:fwd}
\end{figure}

A mathematical description of the transport of a substance (contaminant) concentration $u$ in a bounded open domain $\Omega \subseteq \mathbb{R}^n \text{ for } n\in\{2,3\} $ is given by the following equation:
\begin{equation}\label{eq:forward_equation}
\begin{aligned}
\rpto(u):=  u_t-\diffcoeff\Delta u + \velocity\cdot\nabla u &= 0 &\qquad&\text{in}\ (0,T)\times\Omega,\\
  \diffcoeff\nabla u \cdot \normal &= 0 &&\text{in}\ (0,T)\times (\outflowBoundary \cup \characteristicBoundary),\\
  u&= 0 &&\text{in}\ (0,T)\times \inflowBoundary,\\
  u(0,\cdot) &= \parameterfield &&\text{in}\ \Omega.
\end{aligned}
\end{equation}
The parameter-dependent forward problem shown in~\autoref{eq:forward_equation} is formulated for realizations of the parameter~$\parameterfield$. A visualization of the contaminant dispersion is provided in \autoref{fig:fwd}. 

The underlying transport process is governed by a wind vector field~$\velocity$, which is assumed to be sufficiently smooth, bounded --- i.e., $\velocity \in L^\infty(\Omega, \mathbb{R}^n)$ --- and divergence-free --- i.e., $\nabla \cdot \velocity = 0$. The example wind vector field used hereinafter is shown in \autoref{fig:fwd} (left). Based on the orientation of the wind vector relative to the outward-pointing boundary normal~$\normal$, the domain boundary~$\partial\Omega$ is partitioned into three disjoint subsets: the outflow boundary~$\outflowBoundary \subset \partial\Omega$, where $\velocity \cdot \normal > 0$; the characteristic (or tangential) boundary~$\characteristicBoundary \subset \partial\Omega$, where $\velocity \cdot \normal = 0$; and the inflow boundary~$\inflowBoundary \subset \partial\Omega$, where $\velocity \cdot \normal < 0$, following the convention in~\cite{Elman.2020}.

\subsection{Inverse problem: Source identification} \label{ssec:inverse}
\begin{figure}[hbt!]
\begin{subfigure}{0.48\textwidth}
\includegraphics[width=0.9\linewidth]{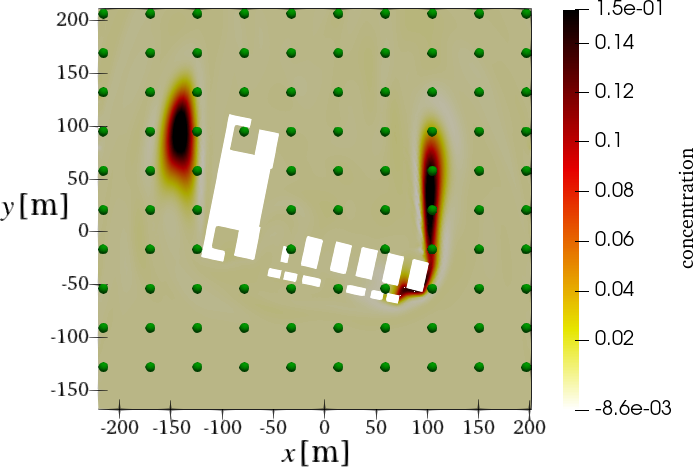}
\end{subfigure}
\begin{subfigure}{0.48\textwidth}
\includegraphics[width=0.9\linewidth]{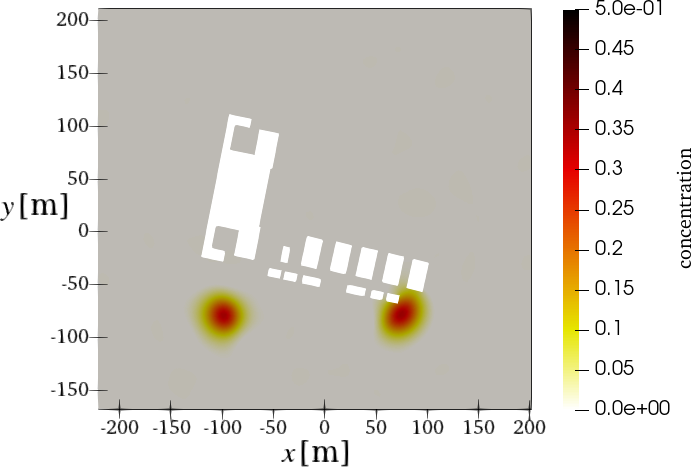}
\end{subfigure}
    \caption{Inverse Problem. Measurements at 96 equidistantly spaced sensor positions (left) and reconstructed initial condition (solution of the inverse problem, right).}
 \label{fig:inverse_2d_pred}
\end{figure}

Whenever measurements of the concentration at discrete locations and times are available, an obvious question is whether the initial condition can be reconstructed on the basis of the given measurements. The respective inverse problem is illustrated in \autoref{fig:inverse_2d_pred}, see also~\cite{Villa.2021}. As function space for the initial condition, we consider a admissible subset of square-integrable functions, for example, $\sop := H_{\inflowBoundary}^{1,2}(\Omega) := \left\{ m \in H^{1,2}(\Omega) \mid m_{|\inflowBoundary} = 0 \right\}$, for this application. In this setting, the estimation of the initial value leads to a linear optimization problem, which will be addressed in the following.
The first step is to describe sensor measurements within this formulation. To do this, we define a well-posed and bounded \textit{space-time observation operator} $\obsO: C^{0}([T_0, T] \times \bar{\Omega}) \to \mathbb{R}^q$ by $u \mapsto \sum_{i=1}^q \delta_{\observationforequation}(u)\,\stdbase_i = \left(u\observationforequation\right)_{i=1}^{q},$ where $\observationforequation \in (T_0, T) \times \Omega$ represents a sequence of space-time coordinates, and $\{\stdbase_i\}$ is the standard basis of $\mathbb{R}^q$. 
Using this observation operator, the final \textit{parameter-to-observable map} $\pto: \sop \to \mathbb{R}^q$ is defined by
\begin{equation}\label{eq:parameter-to-observable}
\pto(\parameterfield) := \obsO \circ \pts(\parameterfield), \quad \text{with } \pts(\parameterfield) := u \text{ such that } \rpto(u) = \parameterfield.
\end{equation}
Here $\pts$ is the \textit{parameter-to-state map}, mapping parameter space $\sop$ to state space, often referred to as 'model' in this context~\cite{Cvetkovic.2024}. 
An example problem in which $\pts$ maps the initial condition $u_0$ to the solution of $u(t=\SI{10}{\second},\cdot)$ is illustrated in \autoref{fig:fwd}. 

The next step is to model sensor noise, which is usually present in real-world measurements. To this end, it is assumed that the observations $\measurement=\mathcal{F}(\parameterfield) + \bm{\epsilon}\,$ contain additive noise $\bm{\epsilon} \thicksim \gauss{\textbf{0}}{\Gamma_{\text{noise}}}$ due to measurement uncertainties. For the sake of simplicity, it is further assumed that the sensor noise at the different sensor positions is uncorrelated and of equal magnitude, represented by the diagonal matrix $\Gamma_{\text{noise}} = \text{diag}(\sigma^2, \dots, \sigma^2)$. Moreover, the conformity of the simulation with the measured values, also called \textit{misfit}, $\by = \mathcal{F}(\parameterfield)-\measurement$, is evaluated in the following norm
\begin{equation*}
\begin{aligned}
\norm{\mathcal{F}(\parameterfield)-\measurement}_{\Gamma^{-1}_\text{noise}}^2 = \fracSolidus{1}{\sigma^2} \sum_{i=1}^{q} \left(u \observationforequation - d_i \right)^2.
\end{aligned}
\end{equation*}
In summary, the inverse problem consists of using available measurements $\measurement$ to infer the values of the unknown parameter field $\parameterfield$. 
Alternatively, the $\mmap$ point can be characterized by the solution of the minimization problem
\begin{equation}\label{eq:objective}
\begin{aligned}
  \mmap=\argmin_{\parameterfield\in \sop} J(\parameterfield):= \frac{1}{2} \norm{\mathcal{F}(\parameterfield)-\measurement}_{\Gamma^{-1}_\text{noise}}^2 + \frac{1}{2}\norm{\parameterfield-\parameterfield_{\text{pr}}}_{\Gamma^{-1}_{\text{prior}}}^2\,,
\end{aligned}
\end{equation}
with the prior information encoded as a Tikhonov regularization term.

Unfortunately, the system is heavily under-determined under real conditions, as sensor measurements are only available at a few locations, but an initial condition is to be reconstructed for the entire domain. In order to transform this into a well-posed problem, some prior knowledge is needed and a formulation as a Bayesian inverse problem provides a suitable framework. 
In this setting, a Gaussian prior $\gauss{\parameterfield_{\text{pr}}}{\covPr}$ with mean $\parameterfield_{\text{pr}}$ and covariance $\covPr$ is chosen for parameter regularization. Then, the posterior density of $\parameterfield$ satisfies by Bayes' theorem:
$
{\pi_{\text{post}}(\parameterfield|\measurement) \propto \pi_{\text{like}}(\measurement|\parameterfield) \, \pi_{\text{prior}}(\parameterfield)}\,$.
Here, $\pi_{\text{like}}(\measurement|\parameterfield) \propto \exp( \frac{1}{2} \norm{\mathcal{F}(\parameterfield)-\measurement}_{\Gamma^{-1}_\text{noise}}^2)$
is the likelihood function under the observational noise  $\bm{\epsilon} \thicksim \gauss{\textbf{0}}{\Gamma_{\text{noise}}}$. Due to the linearity of $\mathcal{F}$, the posterior distribution is again a Gaussian distribution $ \gauss{\mmap}{\Gamma_{\text{post}}}$ with covariance and mean
\begin{equation}
\begin{aligned}
\Gamma_{\text{post}} =  (\mathcal{F}^*\Gamma_{\text{noise}}^{-1}\mathcal{F} + \covPr^{-1})^{-1}  \text{ and } \mmap=\Gamma_{\text{post}}(\mathcal{F}^*\Gamma_{\text{noise}}^{-1}\measurement+ \covPr^{-1}\parameterfield_{\text{pr}})\,.
\end{aligned}\label{eq:mmap}
\end{equation}
The formally adjoint operator \(\mathcal{F}^*: \mathbb{R}^q \to \sop^*\) is required for this framework. For the mapping \(\mathcal{F}: \sop \to \mathbb{R}^q\) between Hilbert spaces, the formal adjoint operator \(\mathcal{F}^*\) is characterized by the relation $
\langle \mathcal{F}(\parameterfield), \by \rangle_{\mathbb{R}^q} = \langle \parameterfield, \mathcal{F}^*(\by) \rangle_{L^2(\Omega)}
$ for all \(\by \in \mathbb{R}^q\) and \(\parameterfield \in \sop\), and its existence follows from Riesz's representation theorem~\cite{Alt.2012}. 
The posterior covariance or Hessian $\hessian:=\pto^*\Gamma_{\text{noise}}^{-1}\pto+ \covPr^{-1}$ of the objective function $J$, see \autoref{eq:objective}, contains a wealth of information about the system. In line with the Bayesian framework~\cite{Villa.2021}, the covariance or inverse Hessian matrix can be employed to predict the uncertainty of the system and is of particular significance for optimal sensor placement in the following chapter.

The mean value $\mmap$ is a reliable estimate for the initial value and thus represents the solution of the inverse problem \autoref{fig:inverse_2d_pred}.
In order to calculate $\mmap$ from \autoref{eq:mmap}, a further specification of the adjoint operator $\pto^*$ is necessary. By applying partial integration to the weak form of \autoref{eq:forward_equation}, the adjoint state $p$
can be derived and satisfies the following equation:
\begin{equation}\label{eq:adjoint_equation}
\begin{aligned}
  -p_t - \diffcoeff \Delta p - \text{div}(p \velocity) &= -\fracSolidus{1}{\sigma^2} \sum_{i=1}^{q}  y_i \, \delta_{\observationforequation} &\qquad & \text{in} \ (0, T) \times \Omega, \\
  (\velocity p + \diffcoeff \nabla p) \cdot \normal &= 0 && \text{in} \ (0, T) \times (\outflowBoundary \cup \characteristicBoundary), \\
  p &= 0 && \text{in} \ (0, T) \times \inflowBoundary, \\
  p(T, \cdot) &= 0 && \text{in} \ \Omega.
\end{aligned}
\end{equation}
for given $\by \in\R^q$. Finally, the adjoint operator $ \pto^*$ can be explicitly determined, resulting in $\pto^*\by=p(0,\cdot)$.

\subsection{Goal-oriented optimal experimental design: Sensor placement} \label{ssec:oed}
So far, a model for the forward problem and an estimate of the initial condition, given a fixed sensor configuration, have been derived, but the question of how such a sensor arrangement ought to be chosen remains unanswered. As mentioned in the previous chapter, the covariance of the posterior, or more specifically, the Hessian matrix, plays a crucial role in developing an indicator for the uncertainty in the system. In classical Bayesian optimal experimental design (OED), the A-optimal design is found throughout the literature. For Gaussian posteriors, the A-optimal design minimizes the trace of the posterior covariance matrix, which amounts to minimizing the average pointwise variance of the inferred parameter
\begin{equation}
\begin{aligned}
\min_{\sensorweights \in \mathbb{W}} \operatorname{tr}[\Gamma_{\text{post}}(\sensorweights)] + \mathcal{R}(\sensorweights)\,,
\end{aligned}\label{eq:sensorobjectiv}
\end{equation}
where $\mathbb{W}$ is the set of all valid sensor configurations and $\mathcal{R}$ is a suitable regularization term. 
As in \cite{Alexanderian.2014}, a finite set of candidate sensor placements $\observationforequation \in \ [T_0,T]\times \Omega$ for $1\leq i\leq q$ will be considered. An example of a spatial grid with $96$ sensors can be seen in \autoref{fig:inverse_2d_pred}. For this set of candidate locations, a weight vector $\sensorweights \in [0,1]^q$ is defined with the $i$-th entry corresponding to the $i$-th location in space and time. In fact, the weight vector decides which measurements are realized or taken into account. In the case of stationary sensors located at positions $\xobs_s$, measurements collected from these spatial points over the entire time horizon, denoted as $({\cdot,\spacepoint_s^\text{obs}})$, are constantly weighted, included in the misfit. 
Hence, the number of independent entries is the weight vector reduces to the number of possible stationary sensor positions. If we consider a mobile sensor, we have a trajectory $\gamma:\{\tobs_0,...,\tobs_s\} \rightarrow \{\xobs_0,...,\xobs_s\}$. For points on the trajectory $(\tobs_i,\xobs_i)$ the sensor weight is $\sensorweights_i=1$. All weights away from the trajectory are set to $0$.
To adjust the forward model to the chosen sensor configuration, we consider the diagonal matrix $\W \in \mathbb{R}^{q \times q}$ with $\W_{ii}=\sensorweights_i$. If we denote the parameter-to-observable map for all sensor positions as $\mathcal{F}$, then for each design $\sensorweights \in \mathbb{W}$, we have
%\begin{equation}\label{eq:Fw}
 $   \mathcal{F}(\sensorweights) = \W \, \mathcal{F}\,.$
%\end{equation}
Taking this further and using $W$ to modify the noise matrix, the influence of a selected sensor layout is also captured in the likelihood function 
\begin{equation*}
    \pi_{\text{like}}(\measurement| \parameterfield, \mathbf{w}) \propto \exp\left\{- \frac{1}{2} (\mathcal{F}(\parameterfield) - \measurement)^T\Wsqrt\Gamma^{-1}_\mathrm{noise} \Wsqrt(\mathcal{F}(\parameterfield) - \bd)\right\}.
\end{equation*}
In consequence, the posterior covariance and mean also depend on the sensor layout via (cf. \autoref{eq:mmap})
\begin{equation*}
\begin{aligned}
\Gamma_{\text{post}}(\sensorweights) =  (\mathcal{F}^*\Wsqrt\Gamma_{\text{noise}}^{-1}\Wsqrt\mathcal{F} + \Gamma_{\text{pr}}^{-1})^{-1}  \text{ and } \mmap(\sensorweights)=\Gamma_{\text{post}}(\mathcal{F}^*\Gamma_{\text{noise}}^{-1}\W \measurement+\Gamma_{\text{pr}}\parameterfield_{\text{pr}})\,.
\end{aligned}
\end{equation*}

To obtain a goal-oriented optimal design, the specialization of the objective function for a quantity of interest (QoI), denoted as $\rho$, must be carried out. For this purpose, another linear operator $\QoiO$ is defined, $\rho=\QoiO\,(\parameterfield)\,,$ which evaluates $\rho$ for a given parameter, here, initial condition. Due to the linearity of $\QoiO$, the prior distribution of $\rho $ is also Gaussian, namely, $\gauss{\rho_{\text{pr}}}{\Sigma_{\text{pr}}}$, with mean $\parameterfield_{\rho}=\QoiO\,(\parameterfield_{\text{pr}})$ and covariance $\Sigma_{\text{pr}}=\QoiO\Gamma_{\text{pr}}\QoiO^*$. This again results in a well-defined Bayesian inverse problem  with a posterior distribution $\pi_{\text{post}}(\rho|\measurement) \thicksim \gauss{\rho_{\text{post}}}{\Sigma_{\text{post}}}$, see~\cite{Attia.2018}. The mean and covariance matrices are given by: $\rho_{\text{post}} =\mathcal{P}\,(\parameterfield_{\text{post}})\,\text{ and } \Sigma_{\text{post}}=\mathcal{P}\,\Gamma_{\text{post}}\,\mathcal{P}^*\,$.

\section{Discretization, preconditioning, and sparsification} \label{sec:method}
\subsection{Finite element discretization}
To solve the PDE problems (\autoref{eq:forward_equation} and \autoref{eq:adjoint_equation}) numerically, a finite element discretization is employed using $\ndof$ Lagrange basis functions $\ansatzSpace = \text{span}\{\phi_1, \dots, \phi_\ndof\}$. Moreover, we find an identification between a vector in $\R^{\ndof}$ and finite elements $\FEMi: \mathbb{R}^\ndof \rightarrow \ansatzSpace$ via $\FEMi(a) = \sum_{i=1}^\ndof a_i \phi_i$. This leads to discretized versions of the parameter-to-observable map $\ptoh:\R^{\ndof}\rightarrow \R^{q}$
defined by \(\ptoh(m_h) = \obsO(u_h)\), where \(u_h\)  solves \autoref{eq:forward_equation} weakly, and its adjoint operator $\ptoh^*:\R^{q} \rightarrow \R^{\ndof}$
given by \(\ptoh^*\mathbf{y} = p_h(0, \cdot)\), where \(p_h\) solves \autoref{eq:adjoint_equation} weakly. The identification is an isometry, i.e., $ \langle \FEMi(a), \FEMi(b) \rangle_{L^2(\Omega)} = \langle a, b \rangle_M=: a^T\,M\,b$, where the mass matrix $M_{ji} := \int_{\Omega} \phi_i(\spacepoint) \phi_j(\spacepoint) \, d\spacepoint,\,M \in \mathbb{R}^{\ndof \times \ndof}$ is used to define the corresponding scalar product. For further details of the finite element discretization, we refer to \cite{Villa.2021} and \cite{Danwitz.02.09.2024}.

\subsection{Preconditioning of discrete inverse problem}
To solve the discrete inverse problem, the prior distribution needs to be carefully chosen. A Laplacian-like operator of trace class $\mathcal{A}:=(\ellipI\,I-\ellipLaplace \Delta)$, with  Robin boundary condition, \linebreak\(\ellipLaplace \nabla m \cdot \normal + \beta m=0 \text{ in }\ (0,T)\times \partial \Omega\), is applied with the constant $\beta$ proposed in~\cite{Daon.2018}. This definition serves as a suitable covariance operator, e.g., $\Gamma_\text{pr}=\mathcal{A}^{-2}=(\ellipI\,I-\ellipLaplace \Delta)^{-2}$. In addition, its discrete counterpart is given by the mapping $\covPrh : \mathbb{R}^{\ndof} \to \mathbb{R}^{\ndof}$ via $\covPrh=(M^{-1}A)^{-2}=A^{-1}MA^{-1}M:=R^{-1}M$,  with matrix representation $A_{ij}=\int\,\phi_i(x)\mathcal{A}\phi_j(x)\,dx$.

Combining this covariance operator together with an appropriate prior mean $m_{\text{pr}}$ (in our applications, e.g., $m_{\text{pr}}=0$) renders the inverse problem well-posed and its solution can be found by solving the following equation for $\mmap$ 
\begin{equation}
\begin{aligned}
\hessianh(\sensorweights)\mmap=\ptoh^*\Gamma_{\text{noise}}^{-1}(W \measurement)+\covPrh^{-1}\,m_{\text{pr}},
\end{aligned}\label{eq:Hmap}
\end{equation}
for the discrete version of the Hessian, that is,  $\hessianh(\sensorweights)=\ptoh^*\Wsqrt\Gamma_{\text{noise}}^{-1}\Wsqrt\ptoh + \covPrh^{-\:1}$. Since determining the Hessian matrix directly is computationally expensive for large-scale problems ($O(\ndof)$-PDE solutions), an iterative conjugate gradient (CG) method is employed. This approach requires only the action of the Hessian-vector on a given vector \( m_k \in \mathbb{R}^{\ndof} \) at each iteration. Specifically, the Hessian action is computed by the following steps: first, solve the forward equation $\measurement = \ptoh(m_k) = \obsO(u_h)$, then, solve the adjoint equation $\ptoh^*\bigl(\Wsqrt \Gamma_{\text{noise}}^{-1} \Wsqrt \, \measurement\bigr) = p_h(0, \cdot)$, next, compute $\tilde{m}_k = \covPrh m_k$, and finally obtain the Hessian action as $\hessianh(\sensorweights) m_k = p_h(0, \cdot) + \tilde{m}_k$.

Since two PDE solutions have to be determined in each iteration, a reduced model of the Hessian matrix is created in advance so that the inverse problem can be solved quickly. 
Using the Cholesky decomposition of the prior covariance, $\covPrh^{-1}=(M^{-1}A)(M^{-1}A)^*$, one obtains the preconditioned Hessian matrix as
\begin{equation}
\begin{aligned}
(\mathcal{A}_h^{-1})^*\hessianh(\sensorweights)\mathcal{A}_h^{-1}=(\ptoh\circ\mathcal{A}_h^{-1})^*\Wsqrt\Gamma_{\text{noise}}^{-1}\Wsqrt(\ptoh\circ\mathcal{A}_h^{-1}) + I,\label{eq:prehessian}
\end{aligned}
\end{equation}
for $\mathcal{A}^{-1}_h=(A^{-1}M)$.
This preconditioned system $\pto\circ\mathcal{A}^{-1}$ has fast decaying eigenvalues and so, we follow \cite{Villa.2021,Isaac.2015} in constructing a low rank approximation of the prior-preconditioned misfit part of the Hessian, i.e., $ \precondhessianhmisfit(\sensorweights):=(\ptoh\circ\mathcal{A}_h^{-1})^*\Wsqrt\Gamma_{\text{noise}}^{-1}\Wsqrt(\ptoh\circ\mathcal{A}_h^{-1})$ by solving the symmetric eigenvalue problem ~\cite{Halko.2011,Liberty.2007}:
\begin{equation*}
\begin{aligned}\hessianhmisfit(\sensorweights)v_i=\lambda_i\,M\,\covPrh^{-1} v_i = \lambda_i Rv_i
\end{aligned}
\end{equation*}
for an orthogonal basis $V_r=(v_1,...,v_r)\in \R^{\ndof \times r}$ and $\lambda_1 \geq ...\geq \lambda_r$ with respect to the scalar product induced by $M\covPrh^{-1}$, i.e., $\langle a, b \rangle_{M\covPrh^{-1}}=: a^T\,M\covPrh^{-1}\,b$. 
Applying the Sherman-Morrison-Woodbury formula, we write
\begin{equation}
\begin{aligned}
\mathcal{A}_h(\hessianh(\sensorweights))^{-1}\mathcal{A}_h&= (\precondhessianhmisfit(\sensorweights)-I)^{-1}\approx I + V_r D_r(\sensorweights) V_r^T
\end{aligned}\label{eq:lowrank}
\end{equation}
where $D_r=\text{diag}(\lambda_1/(1+\lambda_1),...,\lambda_r/(1+\lambda_r))$ is a low rank approximation of the Hessian. Detailed information on this can be found in \cite{Isaac.2015} and \cite{Alexanderian.2014}. Using this approximation, the solution of equation \autoref{eq:Hmap} can be determined with a preconditioned Newton-CG method, see \cite{Villa.2021,Steihaug.1983}.

\subsection{Sparsification of sensor layouts and optimality criteria}
Obviously, the trace of the posterior covariance will be minimal exactly when every sensor weight is set to $1$, which corresponds to using every available piece of information to reduce the level of uncertainty. Thus, to derive a sparse sensor configuration, a penalty term must be introduced into the optimization problem (\autoref{eq:sensorobjectiv}). A common choice for the regularization term in \autoref{eq:sensorobjectiv} is the \(\ell_1\)-norm, which leads to a convex minimization problem with a unique minimizer. Specifically, the regularization term is defined as
\begin{equation}\label{eq:oed_objectiv}
    \mathcal{R}(\mathbf{w}) \coloneqq \alpha \|\mathbf{w}\|_1 = \alpha \mathbf{1}^\top \mathbf{w},
\end{equation}
where \(\alpha > 0\) is the regularization parameter and \(\mathbf{1}\) is a vector of ones. Finally, a binary sensor configuration $\{0, 1\}^{\mathrm{q}}$ is obtained by considering only sensor locations with weights above a certain threshold.

\begin{figure}
\begin{subfigure}{0.48\textwidth}
\includegraphics[width=0.9\linewidth]{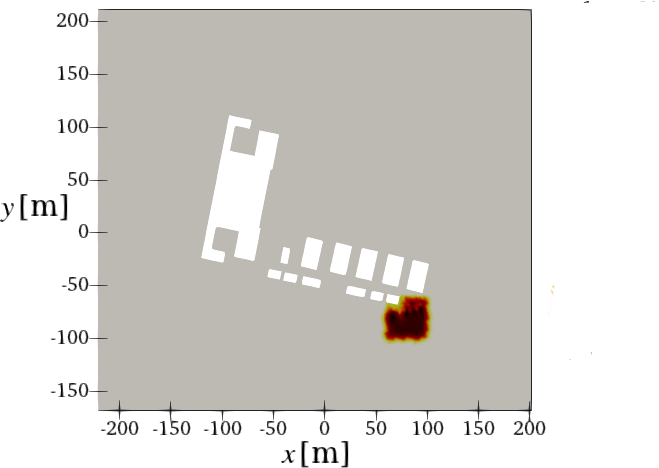}
\end{subfigure}
\begin{subfigure}{0.48\textwidth}
\includegraphics[width=0.9\linewidth]{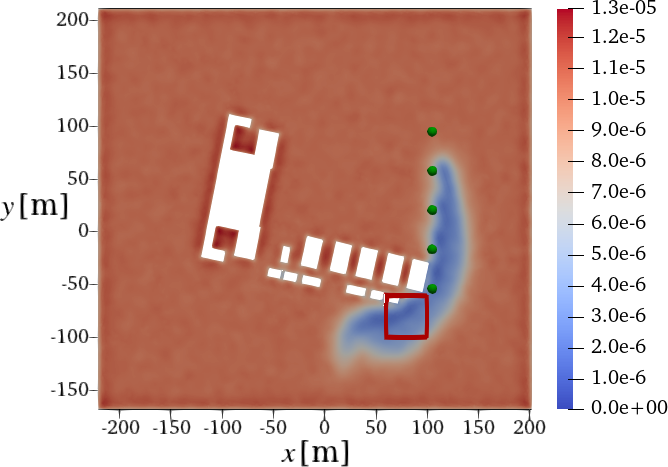}
\end{subfigure}
\caption{Region of interest $\Qoi$ defining the QoI (left) and point-wise variance $\sigma^2 (x)$ of $\mmap$ obtained with five selected sensors (green spheres, right)}
\label{fig:initialQoi}
\end{figure}

The A-optimality criterion for sensor placement minimizes the integrated point-wise posterior variance in linear Bayesian inverse problems: $\int_{\Omega} \mathrm{Var}[m(x)]\,dx = \mathrm{tr}(\Gamma_{\mathrm{post}}(\sensorweights))$. To compute the discrete variance field \(\sigma^2(x) := \mathrm{Var}[m(x)]\), one can extract the diagonal of the inverse Hessian matrix in the finite element basis, assigning each diagonal entry to its corresponding node. However, this exact calculation is computationally expensive (\(O(\ndof)\)), so the reduced-order model (\autoref{eq:lowrank}) is used for visualization. Moreover, the prior covariance is computed using approximate random sampling~\cite{Villa.2021}. Last, but not least, computing the trace of the inverse Hessian is costly. Therefore, we relax A-optimality to a C-optimal design, which only requires evaluating the Hessian’s action on a fixed vector \(c \in \mathbb{R}^{\ndof}\)~\cite{Alexanderian.2021}. \autoref{fig:initialQoi} (left) shows the variance field $\sigma^2 (x)$ of $\mmap$ obtained with a C-optimal sparse sensor layout.

\section{Goal-orientation and sensor steering}
\subsection{Goal-oriented optimal experimental design}\label{ssec:OED}
\begin{figure}
\begin{subfigure}{0.48\textwidth}
\centering
\includegraphics[width=0.9\linewidth]{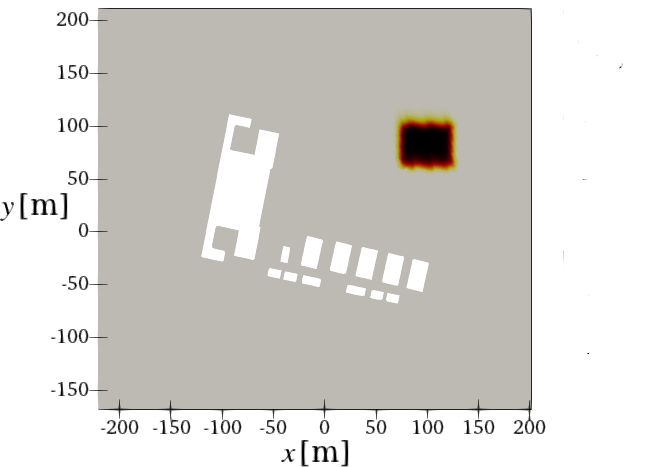}
\end{subfigure}
\begin{subfigure}{0.48\textwidth}
\centering
\includegraphics[width=0.9\linewidth]{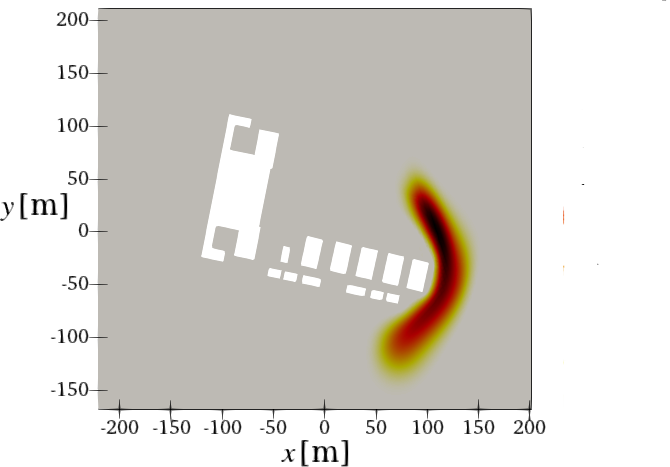}
\end{subfigure}
\caption{Illustration of the QoI $\indFunc{[\Tzeroqoi,\Tfinalqoi]\times \Qoi}$ (left) and solution $c=\pts^*(\indFunc{[\Tzeroqoi,\Tfinalqoi]\times \Qoi})$ of the transport problem induced by the adjoint operator $\pts^*$ (right).}
\label{fig:qoi}
\end{figure}

In the next step, the presented method is adjusted to achieve \emph{goal-oriented} C-optimal experimental designs for stationary sensors. The operator $\QoiO:\sop \rightarrow \R$ is first selected so that the design, namely the sensor placement, is optimized to observe initial conditions in an area $\Qoi \subset \Omega$, shown in \autoref{fig:initialQoi}. This simplification leads to the solution of a C-optimal design. The indicator function of a subset, that is the function $\indFunc{A}\colon X\to \{0,1\}$, which for a given subset $A$ of $\Omega$, attains the value $1$ at points in $A$ and the value $0$ at points outside of $A$. If an optimal design for the initial conditions $m$ is desired, the operator $\QoiO$ does not depend on the solution \(u\), or more specifically, it does not depend on the operator \(\pto\). Concretely, the operator for the QoI is given by $\QoiO(m):=\langle m,\indFunc{\Qoi}\rangle_{L^2(\Omega)}=\int_P m(\spacepoint)\,d\spacepoint$. By identifying the dual space of \(L^2(\Omega)\) and defining the adjoint operator, it follows that \(\QoiO^* = \indFunc{\Qoi}\), and the design function for the optimal experimental design takes the form $ \operatorname{tr}[\Gamma_{\text{post}}(\sensorweights)]= \langle \indFunc{\Qoi}, \Gamma_{\text{post}}(\sensorweights)\indFunc{\Qoi} \rangle_{L^2(\Omega)}.$ In a finite element setting, the function \(\indFunc{\Qoi}\) is represented by a vector \(c \in \mathbb{R}^{n_{\mathrm{dof}}}\) via the usual projection, i.e., \(\indFunc{\Qoi} - c_h \perp \mathcal{V}_h\) into the finite element space \(\mathcal{V}_h = \{\phi_1, \dots, \phi_{n_{\mathrm{dof}}}\}\).
According to \autoref{eq:sensorobjectiv}, the design function is given by $\operatorname{tr}[\Gamma_{\text{post}}(\sensorweights)] =  c_h^T\,M\,\hessianh^{-1}(\sensorweights)\,c_h$.
Together with a suitable regularization term, e.g., $\alpha\lVert \bw \rVert_{l^1}$, an optimal design can be obtained by minimizing the objective $c_h^T\,M\,\hessianh(\sensorweights)\,c_h + \alpha\lVert \bw \rVert_{l^1}$ which is illustrated in \autoref{fig:initialQoi}. If this procedure is generalized to determine an optimal sensor placement for the contaminant concentration over a specific spatial region and time interval, the operator \(\QoiO\) must be extended accordingly. Specifically, by defining the set \(\QoiO\) as a subset of space-time, e.g., \([\Tzeroqoi, \Tfinalqoi] \times \Qoi\) (see \autoref{fig:qoi}), the operator \(\QoiO\) is then constructed as
\begin{equation*}
\begin{aligned}
  \QoiO(\parameterfield)=\int^{0}_{T}\,\int_{\Omega} \indFunc{[\Tzeroqoi, \Tfinalqoi] \times \Qoi}\, \pts(\parameterfield)(t,\spacepoint)\,dt\,d\spacepoint= \int^{\Tfinalqoi}_{\Tzeroqoi}\,\int_{\Qoi} \,u(t,\spacepoint)\,dt\,d\spacepoint.
\end{aligned}
\end{equation*}
Again, the adjoint of $\QoiO$ is needed. For this, we calculate 
\begin{equation*}
\begin{aligned}
\QoiO^*(1)(\hat m) &= \langle \pts(\hat m) , \indFunc{[\Tzeroqoi,\Tfinalqoi]\times \Qoi}\rangle_{L^2(\Omega)} \\
&= \langle \hat m , \pts^*\indFunc{[\Tzeroqoi,\Tfinalqoi]\times \Qoi})\rangle_{L^2(\Omega)},
\end{aligned}
\end{equation*}
and therefore the map $c$ is obtained by $c:=\QoiO^*(1)=\pts^*(\indFunc{[\Tzeroqoi,\Tfinalqoi]\times \Qoi})$, since $\hat m$ was arbitrary. More precisely, the map $c$ satisfies the following PDE:
\begin{equation}\label{eq:adjoint_equation_qoi}
\begin{aligned}
  -c_t - \diffcoeff \Delta c - \text{div}(c \velocity) &= \indFunc{[\Tzeroqoi,\Tfinalqoi]\times \Qoi} &\qquad & \text{in} \ (0, T) \times \Omega, \\
  (\velocity c + \diffcoeff \nabla c) \cdot \normal &= 0 && \text{in} \ (0, T) \times (\outflowBoundary \cup \characteristicBoundary), \\
  c &= 0 && \text{in} \ (0, T) \times \inflowBoundary, \\
  c(T, \cdot) &= 0 && \text{in} \ \Omega.
\end{aligned}
\end{equation}
The solution $c$ is shown for an example case in \autoref{fig:qoi}. 

Returning to the definition, the objective for the goal-oriented sensor design reads
\begin{equation*}
\begin{aligned}
\operatorname{tr}(\Gamma_{\text{post}}(\sensorweights)) &= \QoiO\hessian^{-1}(\sensorweights)\QoiO^*(1) = \QoiO(\hessian^{-1}(\sensorweights) c) \\ &=\langle \pts(\hessian^{-1}(\sensorweights) c),\indFunc{[\Tzeroqoi,\Tfinalqoi]\times \Qoi}\rangle_{L^2(\Omega)} \\
&=\langle \hessian^{-1}(\sensorweights) c , \pts^*(\indFunc{[\Tzeroqoi,\Tfinalqoi]\times \Qoi})\rangle_{L^2(\Omega)}=\langle c, \hessian^{-1}(\sensorweights)c \rangle_{L^2(\Omega)}.
\end{aligned}
\end{equation*}
Thus, the time-dependent case is reduced to the fact that an optimal design for initial conditions for the transported QoI $c$ must be found and coincides with the first case.

For the numerical evaluation of this objective function and its gradient, which are required to minimize \autoref{eq:sensorobjectiv} using the L-BFGS-B solver, we start again with the inverse low-rank approximation of the Hessian from \autoref{eq:lowrank} and proceed to compute the trace as follows
\begin{equation*}
\begin{aligned}
\operatorname{tr}(\Gamma_{\text{post}}(\sensorweights)) &= \langle c, \hessian^{-1}(\sensorweights)c \rangle_{L^2} 
&\approx  c_h^T \mathcal{A}_h^{-1}(I- V_r D_r V_r^T)\mathcal{A}_h^{-1}c_h\,.
\end{aligned}
\end{equation*}
For a shorter notation, we set $\hat q_h := (I- V_r D_r V_r^T)\mathcal{A}_h^{-1}c_h$, $q_h :=\mathcal{A}_h^{-1}\hat q_h$ and obtain $\operatorname{tr}(\Gamma_{\text{post}}(\sensorweights)) \approx c_h^T\, q_h$ and so the calculation of the of the trace consists only a projection in the low rank subspace and solutions of an elliptic problem $\mathcal{A}$ resp. $\mathcal{A}_h$, for which very fast solving strategies exists. To calculate the derivative, we follow~\cite{Alexanderian.2014} and conclude for this simplified case 
\begin{equation}
\begin{aligned}
\ddx{}{\sensorweights_i} \operatorname{tr}(\Gamma_{\text{post}}(\sensorweights))=(\pto^i(q))^2\approx (\ptoh^i\circ\mathcal{A}^{-1}_h(\hat q))^2\label{eq:grad_trace}
\end{aligned}
\end{equation}

So, this calculation can be replaced by a surrogate model for the preconditioned forward operator. In principle, this procedure can be extended to a stronger, A- or D-optimal design.

\subsection{Dynamic sensor steering based on goal-oriented optimal sensor placement}
\label{sec:method_steering}
\begin{figure}
\centering
\begin{subfigure}{0.48\textwidth}
\includegraphics[width=0.9\linewidth]{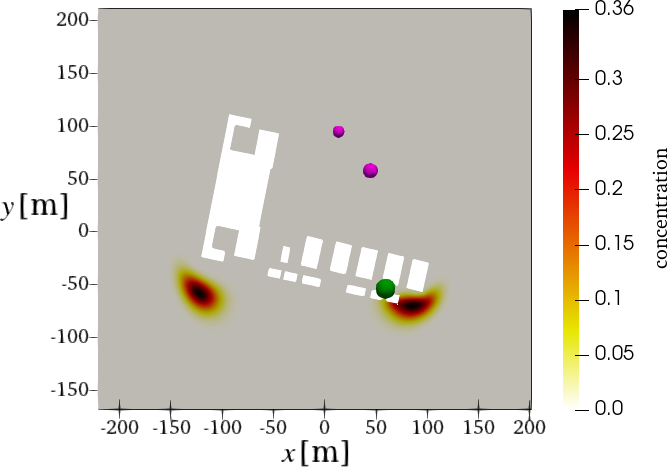}
\caption{Sample information}
\end{subfigure}
\begin{subfigure}{0.48\textwidth}
\includegraphics[width=0.9\linewidth]{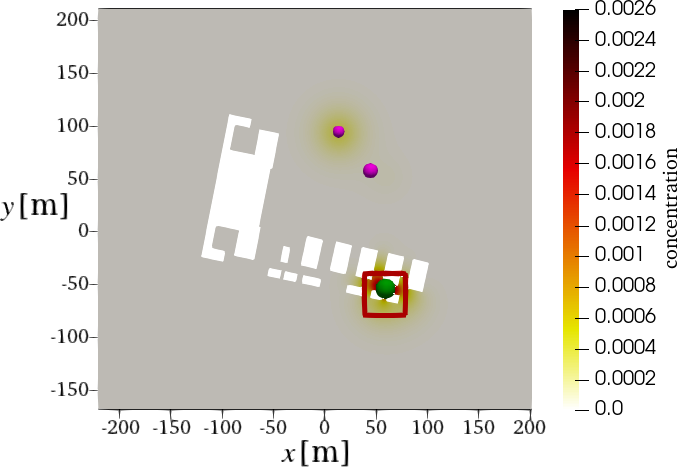}
\caption{Inverse solution}
\end{subfigure}
\begin{subfigure}{0.48\textwidth}
\includegraphics[width=0.9\linewidth]{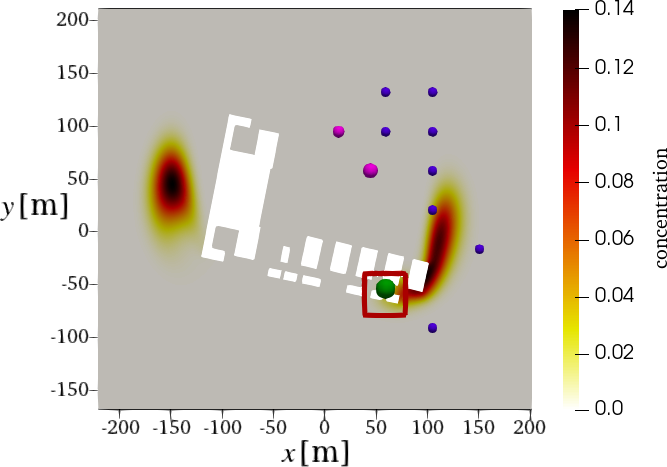}
\caption{Optimal sensor configuration}
\end{subfigure}
\begin{subfigure}{0.48\textwidth}
\includegraphics[width=0.9\linewidth]{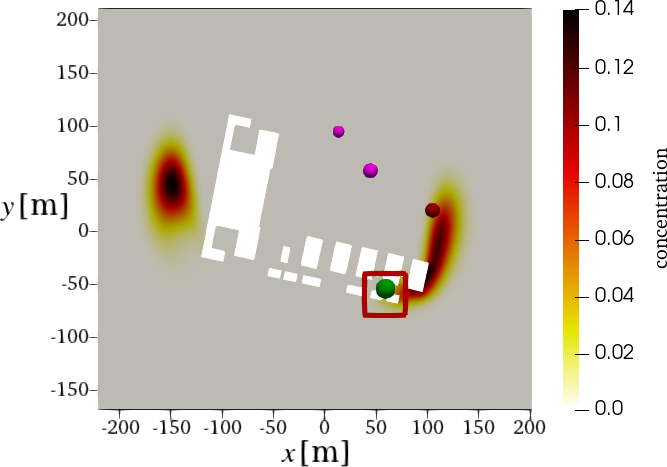}
\caption{Steer sensor}
\end{subfigure}
\caption{Dynamic sensor steering test case with stationary sensor (green) and mobile sensor (purple trajectory). (a) True concentration field at \( u(\cdot, t = \SI{2.2}{\second}) \)   during the second steering time step, (b) Maximum-a-posteriori estimate \( \mmap \) of inverse problem and algorithmically selected zone of interest (red square), (c) optimal sensor design, (d) target position of the mobile sensor (red sphere), true concentration field \( u(\cdot, t = \SI{7}{\second}) \) in background of (c) and (d) to visualize transport problem dynamics}
\label{fig:expl_steering}
\end{figure}

\begin{figure}
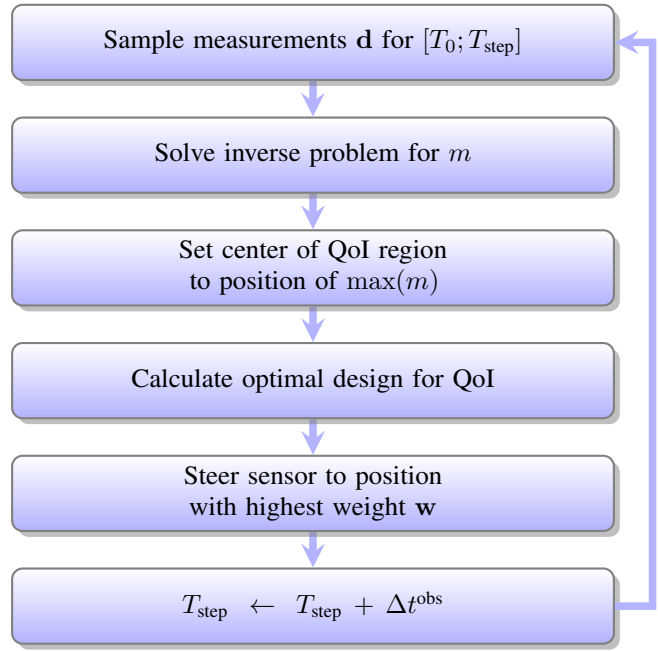

\centering
\smartdiagramset{back arrow disabled=false,uniform color list=blue!30!white for 6 items,
font=\fontsize{10pt}{12pt}\selectfont,
text width=6cm,
module minimum width=8.0cm,
module minimum height= 1.0cm,
module y sep= 1.5,
back arrow disabled= false}

\smartdiagram[flow diagram:vertical]{
Sample measurements $\measurement$ for $[T_0;\Tstep]$,
Solve inverse problem for $\parameterfield$,
Set center of QoI region to position of $\max(\parameterfield)$,
Calculate optimal design for QoI,
Steer sensor to position with highest weight $\sensorweights$,
$\Tstep\leftarrow\Tstep+ \Delta \tobs$
}
\caption{Algorithm for dynamic sensor steering based on goal-oriented sensor placement \cite{wogrin,Wogrin.2023}.}
\label{fig:steeringmethod}
\end{figure}

A method will now be presented which is capable of dynamically controlling a sensor in such a way that a greater knowledge of the actual contaminant concentration can be generated.  
We assume that some knowledge about the concentration is already available due to certain stationary sensors, i.e., that the true contaminant already possesses an appreciable concentration at the sensor location, to permit a solution to the inverse problem. This situation can be seen in \autoref{fig:expl_steering}. We then set the QoI so that its center point is at the maximum of the reconstructed initial condition. The optimum design is then calculated on this basis and the sensor is steered to the position with the highest weight $\sensorweights$. The next measurement is then awaited and the procedure is started again from the beginning. The method is shown schematically in \autoref{fig:steeringmethod}. In this way, we obtain a trajectory $\gamma:\{\tobs_0,...,\tobs_s\} \rightarrow \{\xobs_0,...,\xobs_s\}$ for the steered sensor.

\section{Numerical results of optimal experimental designs} \label{sec:numerical}
In order to simulate scenarios on real-world domains, we use a highly automated process for grid generation. Building imprints as obstacles for two-dimensional contaminant transport are imported directly from Open Street Map (OSM), and locally refined triangular meshes are generated for the region of interest~\cite{Bonari.2024}. The forward model is implemented  with stabilized linear Lagrange finite elements in the software framework \fenics~\cite{Baratta.2023}. The \fenics extension \hippylib (Inverse Problem PYthon library \cite{Villa.2021}) was used in the implementation of the inverse problem. 

In the three following inverse problem OED examples, we use the forward simulation of \autoref{eq:forward_equation} illustrated in \autoref{fig:fwd} as the ground truth. Two radially symmetric functions $$\parameterfield_{\spacepoint_s}(\spacepoint; \spacepoint_s, r) = \min \left\{0.5,\exp \left(-\text{ln}(\epsilon)\norm{\spacepoint-\spacepoint_s}_2^2/r^2\right)\right\}\,,\epsilon=0.001$$
describe the initial concentration field, 
\begin{equation}
 \begin{aligned}
 u_0(\spacepoint)=\parameterfield_{\spacepoint_s}(\spacepoint; \spacepoint_s=[\SI{-100}{m}, \SI{-80}{m}], r = \SI{25}{m})+ \parameterfield_{\spacepoint_s}(\spacepoint; \spacepoint_s= [\SI{75}{m}, \SI{-80}{m}], r= \SI{25}{m}).
 \end{aligned}\label{eq:gauss_blob}
\end{equation}
The initial concentration field is transported by the vector field $\velocity$. For the considered test cases, we estimate a stationary wind vector field as solution of the incompressible Navier-Stockes equations with wind entering the given geometry from south at a velocity of \(\velocity = \SI{10}{\metre\per\second}\). This condition is realized using a Dirichlet boundary condition. In the inner boundaries that represent the imprints of the buildings, a no-slip condition is applied. The remaining edges correspond to free boundary conditions. For the chosen Reynolds number of 50, we obtain the laminar wind field visualized in \autoref{fig:fwd}. Moreover, the diffusion coefficient is selected as $\diffcoeff=\SI{1}{\square \metre \per \second}$ resulting in a transport problem with a moderate Peclet number. Finally, the time step size for the implicit Euler time-stepping scheme is set to $\SI{0.05}{\s}$. In the parametrization of the prior, the constants were chosen as \(\ellipI = 8\) and \(\ellipLaplace = 800\), yielding the operator $\mathcal{A} := 8\,I - 800\,\Delta$.
\begin{figure}
\begin{subfigure}{0.48\textwidth}
\includegraphics[width=0.9\linewidth]{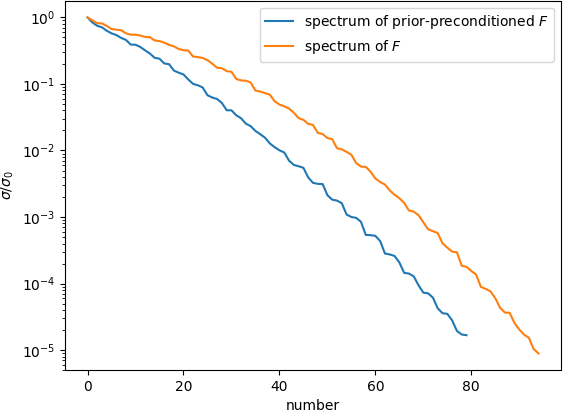}
\end{subfigure}
\begin{subfigure}{0.48\textwidth}
\includegraphics[width=0.9\linewidth]{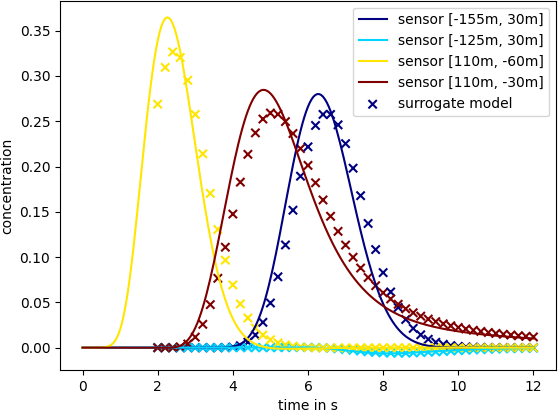}
\end{subfigure}
\caption{Reduced-order modeling. Decay of singular values of $\ptoh$ and of preconditioned $\ptoh\circ\mathcal{A}$ (left), comparison of ROM and forward model $\ptoh$ evaluated at sensor positions (right).}
\label{fig:rom_decay_accuracy}
\end{figure}
In order to make this problem computationally feasible, reduced-order models (ROMs) of the forward and adjoint operators are derived. Considering the forward operator $\ptoh:\R^{\ndof}\rightarrow \R^{q}$, it is observed that it constitutes a linear mapping from a high-dimensional to a lower-dimensional space. Thus, a singular value decomposition is performed to construct a ROM; see also~\cite{Halko.2011,Liberty.2007}. The decomposition provides singular values $\lambda_1 \geq ...\geq \lambda_r$, an $L^2$-orthogonal basis $U_r=(u_1,...,u_r) \in \R^{\ndof \times r}$ and an orthogonal basis $V_r=(v_r,...,v_r) \in \R^{q \times r}$. During the online phase, for example, when the precomputed reduced-order model (ROM) is used for sensor steering, the selected initial condition is projected, such that only matrix-vector multiplications are required:
\begin{equation*}
  \ptoh(m_h) \approx V_r \, D \, U_r \, M \, m_h,
\end{equation*}
where $m_h \in \mathbb{R}^{n_{\mathrm{dof}}}$ and $D = \mathrm{diag}(\lambda_1, \ldots, \lambda_r)$.

As discussed in \autoref{eq:prehessian} and \autoref{eq:grad_trace}, it is also viable to construct a ROM directly for the preconditioned forward operator. The singular values of the operators $\ptoh$ and $\ptoh\circ\mathcal{A}$ are compared in \autoref{fig:rom_decay_accuracy} (left). It is observed that the singular values of the preconditioned operator decay faster and therefore less computations have to be performed to construct a ROM with acceptable accuracy. \autoref{fig:rom_decay_accuracy} (right) shows that the ROM approximates the forward operator fairly well. Furthermore, the computed singular value decomposition is reused to approximate the adjoint operator $\ptoh^*:\R^{q}\rightarrow \R^{\ndof}$ with a ROM as well, namely,$ \ptoh^*(y)\approx M\,U_r\,D\,V_r\,y,$ for $y\in \R^q$.

To assess the benefits of the reduced model in the context of forward evaluations, we first performed a single evaluation of the full-order model, which required approximately $\SI{1}{\second}$ on our hardware. To simulate a scenario relevant to optimal sensor placement, the reduced model was used to evaluate the full-order model at 96 spatial positions over 90 time instances, resulting in a total of $q=8640$ measurements. In total, we calculated 200 spectral values. On average, each evaluation of the reduced model took $\SI{0.00625}{\second}$, yielding a relative speedup of approximately 160 compared to the full-order model.

\subsection{OED 1. Sensor configuration to reconstruct initial condition in critical area}

\begin{figure}
\begin{subfigure}{0.48\textwidth}
\includegraphics[width=0.7\linewidth]{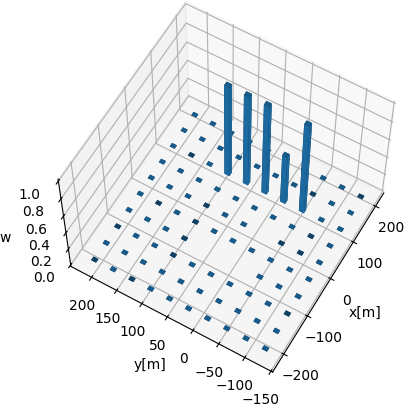}
\end{subfigure}
\begin{subfigure}{0.48\textwidth}
\includegraphics[width=0.9\linewidth]{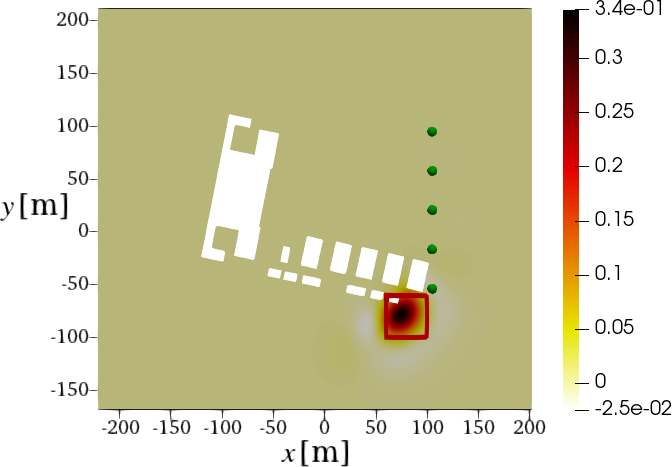}
\end{subfigure}
\caption{OED 1. Weights of optimal sensor configuration to monitor $\Qoi_1$ (left) and reconstructed initial condition (rigth, ground truth shown in \autoref{fig:fwd})}
\label{fig:initial_opt_inverse}
\end{figure}

As first example for a goal-oriented optimal experimental design, we address the problem of identifying an optimal sensor layout to recover the initial condition in a defined subset of the computational domain 
$
\Qoi_1:=\left\{(x,y)\in \Omega\, |\, 75\leq x \leq 125,\, -100 \leq y \leq -60\right\}\,.
$
In a practical application, $\Qoi_1$ might represent a critical area of a chemistry plant site where hazardous material is stored. The inverse problem is posed under the assumption that only stationary sensors are used. These sensors sample the concentration at a rate of \(\SI{5}{\Hz}\), beginning at \(T_0 = \SI{2}{\s}\). Measurements taken after \(\SI{12}{\s}\) are not taken into account. A noise variance of \(\sigma^2 = (0.005)^2\) is assumed, resulting in a signal-to-noise ratio of approximately \(\text{SNR} \approx 100\). Moreover, a regularization parameter of \(\alpha = 0.1\) is applied to obtain a sparse sensor configuration, see \autoref{eq:oed_objectiv}. The selected domain where the quantity of interest (QoI) is inferred is indicated in \autoref{fig:initialQoi} (left), along with the optimal sensor configuration (\autoref{fig:initialQoi} (right)). Moreover, the point-wise variance, which represents the uncertainty in the reconstruction, is also illustrated in \autoref{fig:initialQoi}. The solution to the inverse problem represented by $\mmap$ is visualized in \autoref{fig:initial_opt_inverse}. The numerical result demonstrates a reconstruction quality in $\Qoi_1$ comparable to that achieved using the full configuration with 96 sensors (\autoref{fig:inverse_2d_pred}) in contrast to 5 optimally selected sensors in OED 1.

\subsection{OED 2. Configuration to monitor concentration evolution in critical area}

\begin{figure}
\begin{subfigure}{0.32\textwidth}
\includegraphics[width=0.7\linewidth]{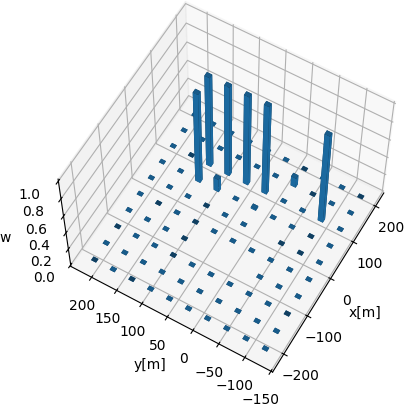}
\end{subfigure}
\begin{subfigure}{0.32\textwidth}
\includegraphics[width=0.9\linewidth]{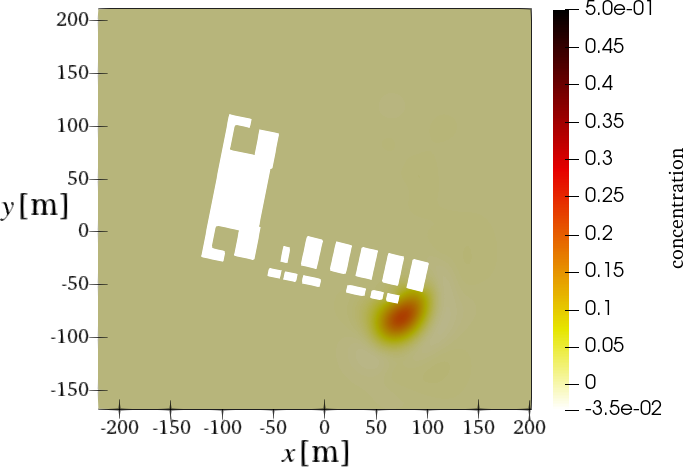}
\end{subfigure}
\begin{subfigure}{0.32\textwidth}
\includegraphics[width=0.9\linewidth]{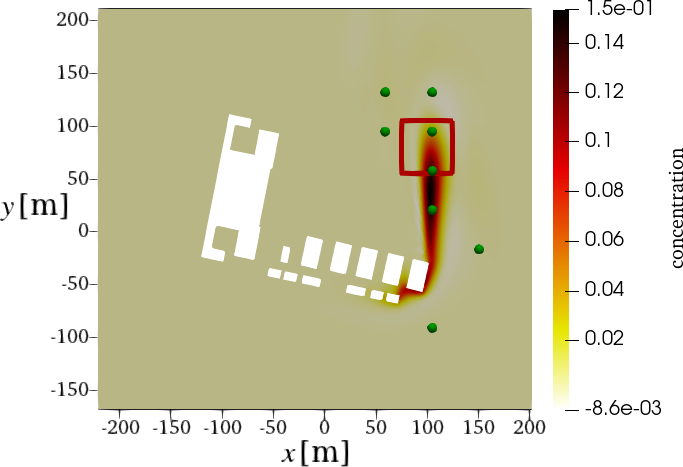}
\end{subfigure}
\caption{OED 2. Weights of optimal configuration to monitor concentration evolution in $\QoiO_2$ (left), reconstructed initial condition (middle), prediction of concentration for $T=\SI{10}{\s}$ (right)}
\label{fig:inverse_opt_time}
\end{figure}

\begin{figure}
\begin{subfigure}{0.48\textwidth}
\includegraphics[width=0.9\linewidth]{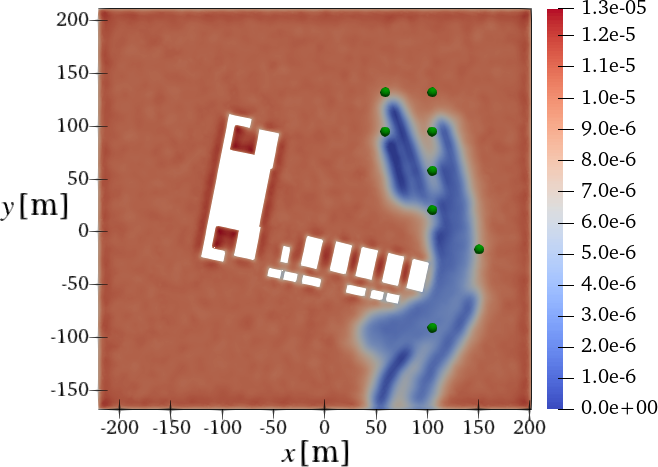}
\end{subfigure}
\begin{subfigure}{0.48\textwidth}
\includegraphics[width=0.9\linewidth]{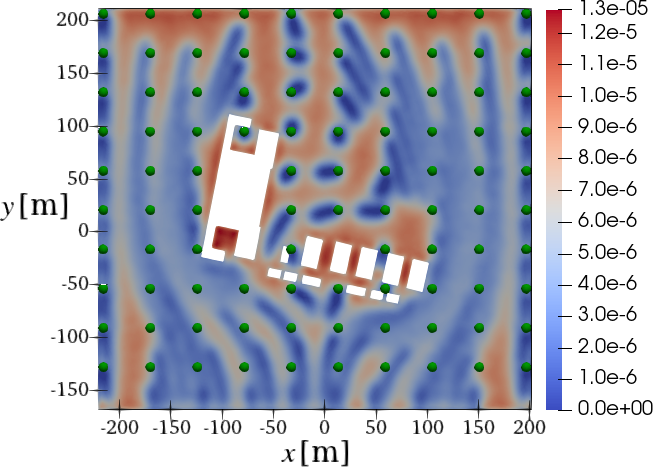}
\end{subfigure}
\caption{OED 2. Point-wise variance $\sigma^2$ as measure of the uncertainty of the inferred parameter for the optimal configuration for $\QoiO_2$ (left) compared to the full sensor grid (right)}
\label{fig:cov_opt_time}
\end{figure}

In the second scenario, our aim is to secure a specific area for a given time period. To achieve this, the quantity of interest (QoI) is defined to depend on the state \(u\). To create a meaningful scenario, we shift $
\Qoi_2:=\left\{(x,y)\in \Omega\, |\, 75\leq x \leq 125,\, 60 \leq y \leq 100\right\}\,.
$ upwards. Goal of OED 2 is to ensure that concentration values can be predicted correctly in region $\Qoi_2$ during the time from $\SI{5}{\s}$ to $\SI{12}{\s}$. This formulation results in the following operator:
 \begin{equation*}
\begin{aligned}
  \QoiO_2(\parameterfield)= \int^{\Tfinalqoi=\SI{12}{\s}}_{\Tzeroqoi=\SI{5}{\s}}\,\int_{\Qoi_2} \,\pts(\parameterfield)(t,\spacepoint)\,dt\,d\spacepoint= \int^{\Tfinalqoi=\SI{12}{\s}}_{\Tzeroqoi=\SI{5}{\s}}\,\int_{\Qoi_2} \,u(t,\spacepoint)\,dt\,d\spacepoint.
\end{aligned}
\end{equation*}
The sensor weights $\sensorweights$ in \autoref{fig:inverse_opt_time} are calculated using a regularization parameter of $\alpha= 1.0$.
In this case as well, the source relevant to the QoI is reconstructed accurately using the optimized sensor configuration. The reconstructed initial condition and the corresponding prediction are shown in \autoref{fig:inverse_opt_time}. As illustrated in \autoref{fig:cov_opt_time}, the reduction in uncertainty is concentrated primarily in the region relevant to QoI compared to the complete sensor configuration.

\subsection{OED 3. Dynamic sensor steering for source identification}

\begin{figure}
\centering
\begin{subfigure}{0.48\textwidth}
\includegraphics[width=0.9\linewidth]{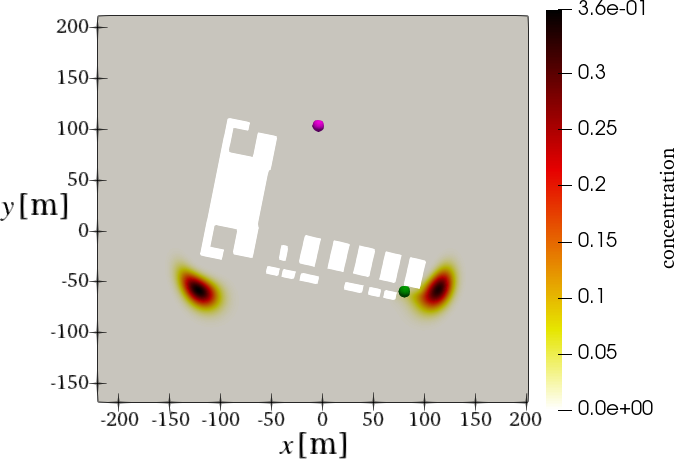}
\caption{Concentration $u$ at $t=\SI{2.0}{\second}$ (truth)}
\end{subfigure}
\begin{subfigure}{0.48\textwidth}
\includegraphics[width=0.9\linewidth]{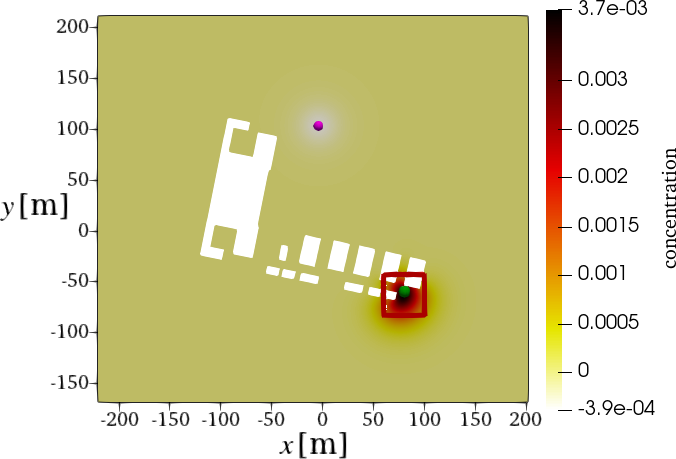}
\caption{Reconstruction of $m$ with data up to $t=\SI{2.0}{\second}$}
\end{subfigure}
\begin{subfigure}{0.48\textwidth}
\includegraphics[width=0.9\linewidth]{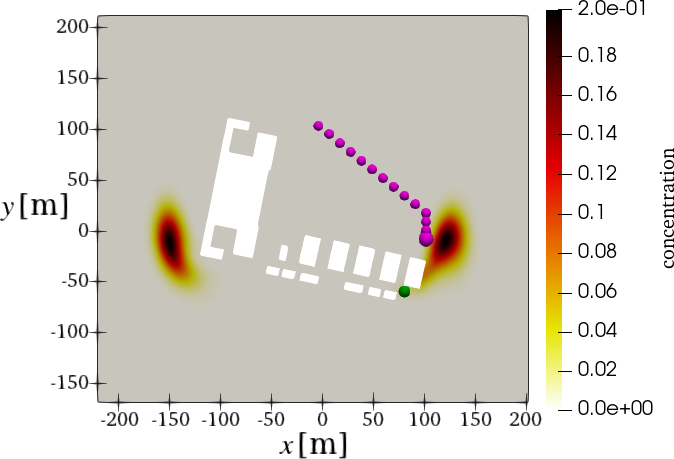}
\caption{Concentration $u$ at $t=\SI{4.6}{\second}$ (truth)}
\end{subfigure}
\begin{subfigure}{0.48\textwidth}
\includegraphics[width=0.9\linewidth]{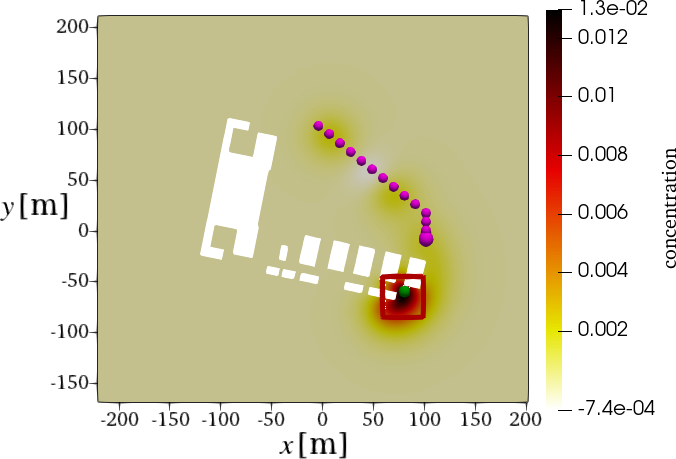}
\caption{Reconstruction of $m$ with data up to $t=\SI{4.6}{\second}$}
\end{subfigure}
\begin{subfigure}{0.48\textwidth}
\includegraphics[width=0.9\linewidth]{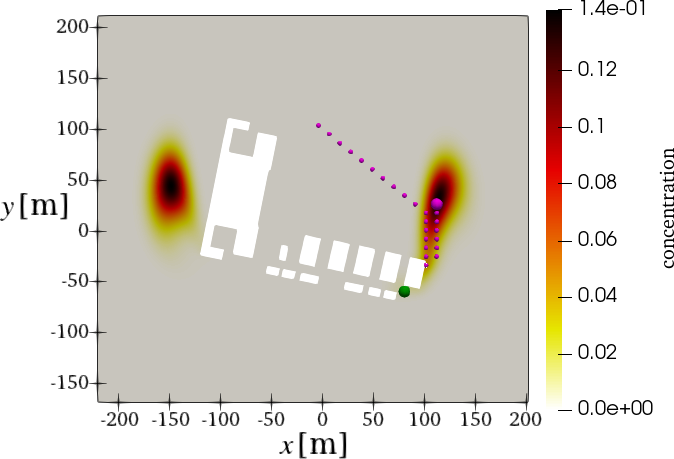}
\caption{Concentration $u$ at $t=\SI{7}{\second}$ (truth)}
\end{subfigure}
\begin{subfigure}{0.48\textwidth}
\includegraphics[width=0.9\linewidth]{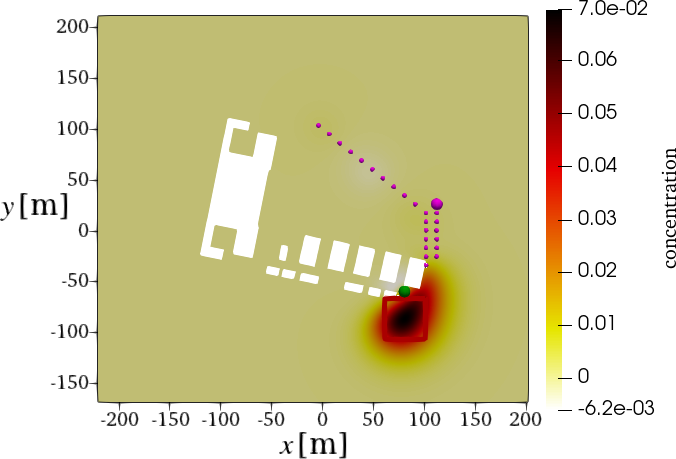}
\caption{Reconstruction of $m$ with data up to $t=\SI{7}{\second}$}
\end{subfigure}
\caption{OED 3. Data fusion of stationary sensor (marked in green) and mobile sensor (marked in purple). Mobile sensor steered to maximize information gain (trajectory marked with small purple spheres)}
\label{fig:steer}
\end{figure}

Finally, the sensor steering method described in \autoref{sec:method_steering} is tested in a numerical application case with $\diffcoeff=\SI{10}{\square \metre \per \second}$. To steer the sensor, a much finer sensor grid totaling 1511 possible sensor locations is used. In this formulation, the moving sensor is allowed to take one step on this grid per cycle, which needs $\SI{0.2}{\second}$ and thus corresponds to the measurement frequency of $\SI{5}{\Hz}$. This corresponds to a speed of approximately $\SI{40}{\metre\per\second}$ for the moving sensor. To demonstrate the capabilities of the sensor steering approach, we placed a single stationary sensor just behind one of the obstacles. However, due to the transport characteristics in this region, information solely from this stationary sensor results in an inaccurate reconstruction of the source, which grossly underestimates the degree of contamination further from the buildings. In addition to the stationary sensor, measurements from a mobile sensor are available. The measurement process begins at time \( T_0 = \SI{2}{\second} \), with data collected at a frequency of $\SI{5}{\Hz}$. The state $u$ at \( \Tstep= T_0 = \SI{2}{\second} \) is shown in \autoref{fig:steer} (a). At this point, the stationary sensor receives very limited information and is thus unable to provide an accurate source estimate. However, computing the optimal sensor design based on the current quantity of interest (QoI), defined as the integral over a square measuring $\SI{40}{\metre}$ on each side, centered on the maximal point of the reconstructed initial condition, yields favorable estimates for informative measurement positions. The mobile sensor is then directed toward the location associated with the highest weight \( \sensorweights \), as determined by the C-optimal design criterion, computed over the observation period $[\Tstep, \Tstep + \SI{2}{\second}]$ with the same sampling rate of $\SI{5}{\Hz}$, wherein we take as sensor weights in the time leading up to \(\Tstep\), the actual past locations of the sensor.
In the subsequent time steps, illustrated in~\autoref{fig:steer} (b), the sensor continues to move toward regions of increasing concentration. In~\autoref{fig:steer}(d), corresponding to \( \Tstep = \SI{4.6}{\second} \), the mobile sensor has found the core of the contaminant.
Finally,~\autoref{fig:steer}(e) and (f) demonstrate that the mobile sensor continues to accurately pursue the contaminant in further time steps. 

\begin{figure}
\centering
\begin{subfigure}{0.48\textwidth}
\includegraphics[width=0.9\linewidth]{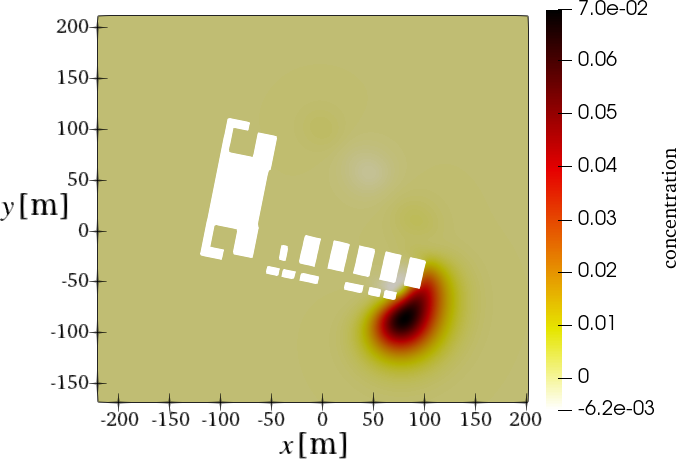}
\caption{Reconstruction of $m$ (with mobile sensor)}
\end{subfigure}
\begin{subfigure}{0.48\textwidth}
\includegraphics[width=0.9\linewidth]{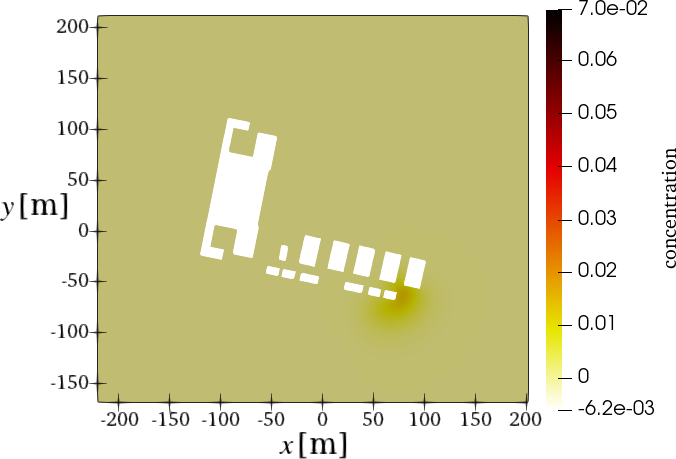}
\caption{Reconstruction of $m$ (stationary sensor alone)}
\end{subfigure}
\begin{subfigure}{0.48\textwidth}
\includegraphics[width=0.9\linewidth]{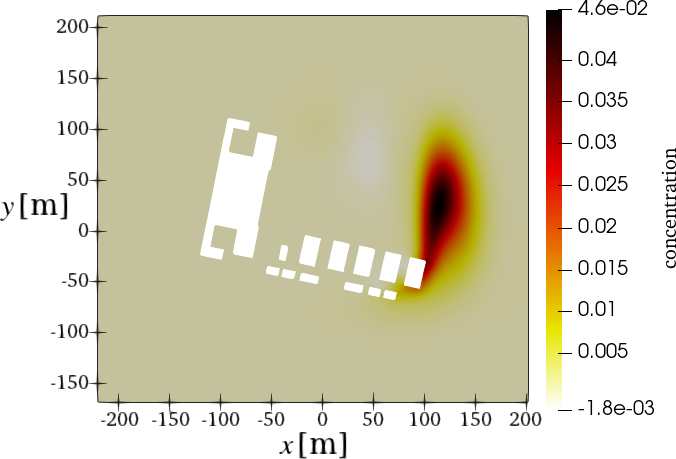}
\caption{Prediction $t=\SI{7}{\second}$ (with mobile sensor)}
\end{subfigure}
\begin{subfigure}{0.48\textwidth}
\includegraphics[width=0.9\linewidth]{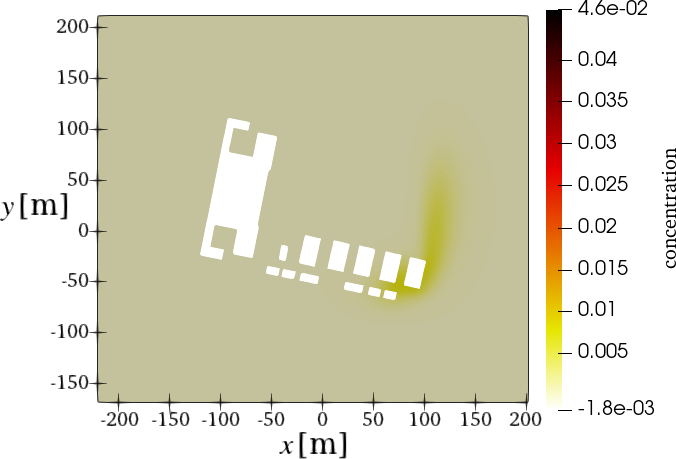}
\caption{Prediction  $t=\SI{7}{\second}$ (stationary sensor alone)}
\end{subfigure}
\caption{OED 3. Comparison of predictions based on data fusion of mobile sensor and stationary sensor (left column) and stationary sensor alone (right column)}
\label{fig:steerCompare}
\end{figure}
Comparing the performance of the stationary sensor with that of the combination of a stationary sensor and a dynamically steered one, we find, as depicted in~\autoref{fig:steerCompare}, that the mobile sensor achieves substantially improved reconstruction accuracy after just \( \SI{7}{\second} \), in contrast to the stationary sensor, which fails to produce a reliable estimate even after \( \SI{12}{\second}\).

\section{Conclusion and Outlook}\label{sec:conc} 

This paper investigated a novel approach for goal-oriented optimal static sensor placement and dynamical sensor steering for PDE-constrained problems. Adopting previous work by Wogrin et al.~\cite{Wogrin.2023,wogrin} on dynamic sensor steering, we leverage a Bayesian approach for the solution of the inverse problem, accelerated by offline-computations of low-rank approximations for the Hessian matrix and an online preconditioned inexact Newton-CG method. The resulting framework was then applied to a more complex geometry extracted from real-world map data.
We showcased the strengths of the proposed workflow on three application cases from the field of airborne contaminant transport: In the first example, we derived an optimal placement of stationary sensors to recover the initial condition inside a spatially constrained rectangular region. The results showed that our proposed method only requires five sensors to reconstruct the initial condition locally with an accuracy comparable to the full configuration of 96 sensors. In our second example, we extended the QoI in the sense that a region of interest is monitored not only at a specific time instance, but over a fixed time period. From a practical point of view, this corresponds to the goal of securing a specific area for a given time period. Using only eight sensors, the evolution of the concentration was accurately reconstructed and the uncertainty was minimized in the area of interest. Lastly, we investigated a dynamic sensor steering problem. Here we showed, that while we are still able to roughly predict the general shape of the initial condition with a unfavorably placed stationary sensor, adding a mobile sensor we obtain much better agreement with the true solution while simultaneously reducing the required measurement time to one-third of the stationary case. This proves that the presented method is able to successfully handle the complexity of a moving sensor and steer the sensor to achieve a fast and reliable reconstruction of the (in practice unknown) initial condition. 

While we believe this work to be an important step towards optimally steering unmanned sensor platforms in crisis situations, there still remain several points for improvement and further investigation of the proposed algorithm. One point for improvement lies within the solution of the inverse problem. One can employ the reasonable assumption that the initial condition is typically sparse in the considered applications. Integrating this additional knowledge into the solution procedure is expected to speed up the time to solution and further improve real-time capabilities the method~\cite{Pieper.2021}. Moreover, we plan to extend the sensor steering to Reinforcement Learning based approach, where the position and size of the QoI in each step is determined by an agent that was previously trained based on trial-and-error interactions with the forward model~\cite{Fricke.2023}.

\section*{Acknowledgement}

DW and AP gratefully acknowledge the funding by dtec.bw - Digitalization and Technology Research Center of the Bundeswehr (project RISK.twin). dtec.bw is funded by the European Union - NextGenerationEU.

\bibliographystyle{unsrtnat}
\bibliography{references}

\end{document}